\documentclass{amsart} 
\usepackage{amssymb} 
\usepackage{enumerate}

\def\np{\medskip} 
\def\nl{\smallskip\noindent}

\def\Sym{{\rm S}} 
\def\Supp{{\rm Supp}} 
\def\Geod{{\rm Geod}} 
\def\Proj{{\rm Proj}} 
\def\M{{\rm M}} 
\def\wh{\widehat } 

\newtheorem{Thm}{Theorem}[section] 
\newtheorem{Lm}[Thm]{Lemma} 
 
\newtheorem{Def}{Definition}

\newtheorem{Prop}[Thm]{Proposition} 
\newtheorem{Cor}[Thm]{Corollary} 

\theoremstyle{definition} 
\newtheorem{defn}[Thm]{Definition} 
 
\newtheorem{Remarks}[Thm]{Remarks}

\newtheorem{Alg} [Thm]{Algorithm}

\renewcommand{\phi}{\varphi} 
\DeclareMathOperator{\tr}{tr} 
\DeclareMathOperator{\het}{ht} 
 
\def\Q{{\mathbb Q}} 
\def\R{{\mathbb R}} 

\newcommand{\ident}{1}

\newcommand{\C}{{\mathbb{C}}} 
 
\newcommand{\N}{{\mathbb{N}}}

\newcommand{\A}{{\rm A}} 
\newcommand{\D}{{\rm D}} 
\newcommand{\E}{{\rm E}} 
\newcommand{\ADE}{{\rm ADE}} 
\newcommand{\adj}{\sim}

\def\a{\alpha} 
\def\b{\beta} 
\def\c{\gamma} 
\def\s{\sigma}

\author{Arjeh M. Cohen 
\& Di\'{e} A.H. Gijsbers
\& David B. Wales} 
\address{Arjeh M. Cohen\\ 
Department of Mathematics and Computer Science\\ 
Eindhoven University of Technology\\ 
POBox 513\\ 
5600 MB Eindhoven\\ 
The Netherlands} 
\email{A.M.Cohen@tue.nl} 
\address{Di\'{e} A.H. Gijsbers\\
Department of Mathematics and Computer Science\\ 
Eindhoven University of Technology\\ 
POBox 513\\ 
5600 MB Eindhoven\\ 
The Netherlands} 
\email{D.A.H.Gijsbers@tue.nl}
\address{David B. Wales\\ 
Mathematics Department\\ 
Sloan Lab\\ 
Caltech\\ 
Pasadena, CA 91125\\ 
USA} 
\email{dbw@its.caltech.edu} 

\title{BMW algebras of simply laced type} 

\begin{document} 

\begin{abstract}
  It is known that the recently discovered representations of the Artin groups
  of type $\A_n$, the braid groups, can be constructed via BMW algebras. We
  introduce similar algebras of type $\D_n$ and $\E_n$ which also lead to the
  newly found faithful representations of the Artin groups of the
  corresponding types. We establish finite dimensionality of these algebras.
  Moreover, they have ideals $I_1$ and $I_2$ with $I_2\subset I_1$ such that
  the quotient with respect to $I_1$ is the Hecke algebra and $I_1/I_2$ is a
  module for the corresponding Artin group generalizing the Lawrence-Krammer
  representation.  Finally we give conjectures on the structure, the dimension
  and parabolic subalgebras of the BMW algebra, as well as on a generalization
  of deformations to Brauer algebras for simply laced spherical type other than $\A_n$.
\end{abstract} 

\maketitle

\section{Introduction} 
In \cite{CW}, representations were given for the Artin groups of spherical
type which are faithful, following the construction of Krammer for braid
groups, \cite{Krammer}. 
(We note that \cite{bigelow} also contains a proof of
the faithfulness of this representation for type $\A_n$,
and that \cite{digne} also generalizes this result to arbitrary spherical types.) 
 Faithful representations for the Artin groups of type
$\A_n$, $\D_n$, and $\E_m$ for $m=6,7,8$ were explicitly constructed. Since
each Artin group of spherical irreducible type embeds into at least one of these,
this shows each is linear. As the representations for type $\A_n$ occur in
earlier work of Lawrence \cite{lawrence}, they are called Lawrence-Krammer
representations.

Zinno, \cite{Zinno}, observed that the Lawrence-Krammer representation of the
Artin group of type $\A_n$, the braid groups on $n+1$ braids, factors through
the BMW algebra, the Birman-Murakami-Wenzl algebra introduced in
\cite{wenzl,murakami}.

In this paper we introduce algebras similar to the BMW algebra for other
types. We associate a unique algebra with each simply laced Coxeter diagram
$M$ of rank $n$.  Here, simply-laced means that $M$ has no multiple bonds. We
define the algebras by means of $2n$ generators and five kinds of relations.
For each node
$i$ of the diagram $M$ we define two generators $g_i$ and $e_i$
with $i = 1,\ldots,n$. If two nodes are connected in the diagram we write $i
\sim j$, with $i,j$ the indices of the two nodes, and if they are not
connected we write $i \not \sim j$.  We let $l,x$ be two indeterminates.

\begin{Def}\label{BMW-def}
Let $M$ be a simply laced Coxeter diagram of rank $n$. The BMW algebra of
type $M$ is the algebra, denoted by $B(M)$ or just $B$, with unit element, over $\Q(l,x)$, 
whose presentation is given on generators $g_i$ and
$e_i$ ($i=1,\ldots,n$) by the following defining relations

\begin{center}
\begin{tabular}{lcllr}
(B1)&\qquad& $g_ig_j=g_jg_i$& when &$i \not \sim j$,\\
(B2)&\qquad& $g_ig_jg_i=g_jg_ig_j$& when &$i\sim j$,\\
(D1)&\qquad& $me_i=l(g_i^2+mg_i-1)$& for all &$i$,\\
(R1)&\qquad& $g_ie_i=l^{-1}e_i$& for all &$i$,\\
(R2)&\qquad& $e_ig_je_i=le_i$& when &$i \sim j$,
\end{tabular}
\end{center}
where $m = (l-l^{-1})/(1-x)$.
\end{Def}

\np Notice that the first two relations are the braid relations commonly
associated with the Coxeter diagram $M$. 
Just as for
Artin and Coxeter groups, if $M$ is the disjoint union of two diagrams $M_1$
and $M_2$, then $B$ is the direct sum of the two BMW algebras $B(M_1)$ and
$B(M_2)$. For the solution of many problems concerning $B$, this gives an easy reduction to
the case of connected diagrams $M$.

In (D1) the generators $e_i$ are expressed in terms of the $g_i$ and so $B$ is
in fact already generated by $g_1,\ldots,g_n$. We shall show below that the
$g_i$ are invertible elements in $B$, so that there is a group homomorphism
from the Artin group $A$ of type $M$ to the group $B^\times$ of invertible
elements of $B$ sending the $i$-th generator $s_i$ of $A$ to $g_i$. As we
shall see at the end of Section \ref{sec:krammer}, the Lawrence-Krammer
representation is a constituent of the regular representation of $B$. 
This generalizes Zinno's result \cite{Zinno}. As
a consequence of \cite{CW}, the homomorphism $A\to B^\times$ is injective.

The fact that the BMW algebras of type $\A_n$ coincide with those defined by
Birman \& Wenzl \cite{wenzl} and Murakami \cite{murakami} is given in Theorem
\ref{wenzl-thm}.

The Lawrence-Krammer representation of the Artin groups is based on two
parameters, in \cite{CW} denoted by $t$ and $r$. The two parameters $m$ and
$l$ here are related by $m=r-r^{-1}$ and $l=1/(tr^3)$.

Our first major result is as follows.

\begin{Thm}\label{findim-thm}
The BMW algebras of simply laced spherical type
are finite dimensional.
\end{Thm} 

The proof is at the end of \S\ref{sec:prelims}.  Some information and 
conjectures about dimensions appear in Section \ref{sec:concon}.

Let $I_1$ be the ideal of $B$ generated by all $e_i$, and let $I_2$ be the
ideal generated by all products $e_ie_j$ for $i$ and $j$ distinct and not
connected in $M$. Then clearly $I_2\subseteq I_1$. Moreover,
it is immediate from the defining relations of $B$ that
$B/I_1$ is the Hecke algebra of type $M$. 
The main result of this paper concerns the structure of $I_1/I_2$.

Let $(W,R)$ be the Coxeter system of type $M$.  We write $\Phi^+$ for the set
of positive roots of the Coxeter system of type $M$.  By $\a_0$ we denote its
highest root, and by $C$ the set of nodes $j$ in $M$ with $(\a_j,\a_0) = 0$
(so $C$ has corank $2$ if $M$ is of type $A$ and $1$ otherwise).  If $X$ is a
set of nodes of $M$, we denote by $W_X$ the parabolic subgroup of $W$
corresponding to $X$. This means that $W_X$ is the subgroup of $W$ generated
by all $r_j$ for $j\in X$.

\begin{Thm} \label{main-thm}
Let $B$ be the BMW algebra of type $\A_n$ $(n\ge1)$, $\D_n$ $(n\ge4)$, or $\E_n$ $(n=6,7,8)$.
Then $B/I_2$ is semi-simple over $\Q(l,x)$.
Let $Z_0$ be the Hecke algebra of type $C$.
For each irreducible representation $\theta$ of $Z_0$,
there is a corresponding representation $\Gamma_\theta$ of $B$ of dimension
$ |\Phi^+|\, \dim(\theta)$ and, up to equivalence, these are the
irreducible representations of $B$ occurring in $I_1/I_2$. 
In particular, the dimension of $I_1/I_2$ as a vector space over
$\Q(l,x)$ equals $|\Phi^+|^2\,|W_C|$.
\end{Thm}

\np The proof of the theorem consists of two major parts.  In
Section \ref{sec:I1modI2}, we provide, for each node $i$ of $M$, a linear spanning
set for $I_1/I_2$ parametrized by triples consisting of two positive roots and an element
of $W_C$.  This shows that $|\Phi^+|^2 |W_C|$ is an upper bound for the
dimension of $\dim(I_1/I_2)$.  The proof that the same number is a lower bound
takes place in Section \ref{sec:krammer}, where the Lawrence-Krammer representation
of $A$, studied in \cite{CW}, is generalized to a representation of the same
dimension as before, viz.~$|\Phi^+|$, but now over the non-commutative ring of
scalars $Z_0$.  Up to a field extension of the scalars, $Z_0$ is well known to
be isomorphic to the group algebra of $W_C$, so $\dim(Z_0) = |W_C|$.

In the final section, we discuss how the results might carry over to $I_2$ and
for $I_r$ with $r\geq 3$. We give a conjecture for the dimension of the BMW algebras of
types $\D_n$ $(n\ge4)$ and $\E_n$ $(n=6,7,8)$.  In the theory of Coxeter
groups and Artin groups, there is a notion of standard parabolic subgroups.
These are subgroups generated by a subset $J$ of the nodes of $M$ and have the
special property that they are Coxeter, respectively, Artin groups of type
$M|_{J}$.  We expect that, at least for spherical $M$, the subalgebra of $B$
generated by the $g_j$ for $j\in J$ will be isomorphic to the BMW algebra of
type $M|_{J}$. For type $\A_n$, the Brauer algebra, cf.~\cite{brauer}, is
obtained as a deformation of the BMW algebra. We conjecture that a similar
deformation exists for the spherical simply laced types, in which the
`pictures', forming the monomial basis of the Brauer algebra, are indexed by a
combinatorial generalization of the abovementioned triples.  As a consequence
of Theorem \ref{main-thm}, these conjectures hold for the quotient algebra
$B/I_2$.  We also discuss possible extensions to other spherical types.

The properties of Artin groups needed for the study of our algebras, are
mentioned in Section \ref{sec:artin}.  The subsequent section
contains a discussion of ideals.
We begin however by studying direct consequences of the defining relations.

\section{Preliminaries}\label{sec:prelims}
For the duration of this section, we let $M$ be a simply laced Coxeter diagram
of rank $n$, and we let $B$ be the BMW algebra of type $M$ over $\Q(l,x)$.

The following proposition collects several identities that are useful for
the proof of the finite dimensionality of $B$, Theorem \ref{findim-thm}. 
Recall that $m$ is related to $x$ and $l$ via
\begin{eqnarray}\label{def-x-eq}
m&=& (l-l^{-1})/(1-x).
\end{eqnarray} 

\begin{Prop}\label{ids-i-prop}
For each node $i$ of $M$, the element $g_i$ is invertible in $B$ 
and the following identities hold.
\begin{eqnarray}
\label{eigi-eq}
e_ig_i &=& l^{-1}e_i,\\
\label{g-inv-eq}
g_i^{-1}&=&g_i+m-me_i,\\
\label{g-sq-eq}
g_i^2&=&1-mg_i+ml^{-1}e_i,\\
\label{e-sq-eq}
e_i^2 &=& xe_i.
\end{eqnarray} 
\end{Prop}

\begin{proof}
By (D1), $e_i$ is a polynomial in $g_i$, so $g_i$ and $e_i$ commute,
so (\ref{eigi-eq}) is equivalent to (R1).

{From} (D1) we obtain the expression $g_i^2 + mg_i -ml^{-1}e_i = 1$.
Application of (R1) to the third monomial on the left hand side gives $g_i(g_i+m-me_i)
= 1$.
So $g_i^{-1}$ exists and is equal to $g_i+m-me_i$.
This establishes (\ref{g-inv-eq}).

Also by (D1), the element $g_i^2$ can be rewritten to a linear combination of $g_i$, $e_i$
and $1$, which leads to (\ref{g-sq-eq}).

As for (\ref{e-sq-eq}),
using (D1) and (R1), we find
$$
e_i^2=e_i \, lm^{-1}(g_i^2+mg_i-1)
= lm^{-1}(l^{-2}e_i+ml^{-1}e_i-e_i) 
=xe_i. 
$$
\end{proof}

\begin{Remarks} 
\label{opLm}
(i). There is an anti-involution on $B$ determined by
$$g_{i_1}\cdots g_{i_q} \mapsto g_{i_q}\cdots g_{i_1}$$
on products of generators $g_i$ of $B$.
We denote this anti-involution by $x\mapsto x^{op}$.

\nl(ii). The inverse of $g_i$ can be used for a different definition of the $e_i$, namely
$$ e_i= 1 + m^{-1}(g_i-g_i^{-1}) \mbox{ \rm for all }i.$$ 

\nl(iii).
By (\ref{e-sq-eq}), the element $x^{-1}e_i$ is an idempotent of $B$ for each
node $i$ of $M$.
\end{Remarks}

The braid relation (B2) for $i$ and $j$
adjacent nodes of $M$ can be seen as a way to rewrite an
occurrence $iji$ of indices into $jij$. It turns out that there are more
of these relations in the algebra, with some $e$'s involved. 

\begin{Prop}\label{iji=jij-id-lm}
The following identities hold for $i\sim j$.
\begin{eqnarray}
\label{gjgiej-eq} \label{lm1}
g_jg_ie_j&=&e_ig_jg_i=e_ie_j,\\
\label{gjgiej-inv-eq} \label{gjeigj-eq} 
g_je_ig_j&=&g_i^{-1}e_jg_i^{-1}\\
\nonumber{} &=& g_ie_jg_i+ m(e_jg_i - e_ig_j + g_ie_j - g_je_i) + m^2(e_j-e_i) \\
\label{8-eq}\label{ejeigj-eq-eq}
e_je_ig_j&=&e_jg_i^{-1} =  e_jg_i + m(e_j-e_je_i),\\
\label{7-eq}
g_je_ie_j&=&g_i^{-1}e_j =  g_ie_j + m(e_j-e_ie_j). \\
\label{eiejei-eq}\label{B4}
e_ie_je_i&=&e_i.
\end{eqnarray} 
\end{Prop}

\begin{proof}
By (D1) and (B2),
\begin{eqnarray*}
g_jg_ie_j & = & g_jg_i(lm^{-1}(g_j^2+mg_j-1)) 
 =  lm^{-1}(g_ig_jg_ig_j+mg_ig_jg_i-g_jg_i) \\
& = & lm^{-1}(g_i^2g_jg_i+mg_ig_jg_i-g_jg_i) 
 =  lm^{-1}(g_i^2+mg_i-1)g_jg_i \\
& = & e_ig_jg_i ,
\end{eqnarray*}
proving the first equality in (\ref{gjgiej-eq}).

We next prove
\begin{eqnarray}\label{lm2}
e_ig_j^ng_ie_je_i &=& e_ig_j^{n-1}e_i
\mbox{ for }n \in \N,\ n\ge1.
\end{eqnarray} 
Indeed, by (B2), (R1), (R2), and the first identity of
(\ref{gjgiej-eq}),
which we have just established, 
$$
e_ig_j^ng_ie_je_i  =  e_ig_j^{n-1}(e_ig_jg_i)e_i 
 =  e_ig_j^{n-1}e_ig_j(g_ie_i) 
 =  l^{-1}e_ig_j^{n-1}e_ig_je_i 
 =  e_ig_j^{n-1}e_i .
$$
The following relation is very useful for
determining relations between the $e_i$.
\begin{eqnarray} \label{lm4}
e_ie_jg_ie_je_i &=& (l+m^{-1})e_i - m^{-1}e_ie_je_i.
\end{eqnarray}
To verify it, we start rewriting one factor $e_j$ by means of (D1),
and then use (\ref{lm2}) with $n=2$ and $n=1$ as well as (R1) and (R2):
\begin{eqnarray*}
e_ie_jg_ie_je_i & = & e_i(lm^{-1}(g_j^2+mg_j-1))g_ie_je_i 
 =  lm^{-1}(le_i + mxe_i - l^{-1}e_ie_je_i) \\
& = & (l+m^{-1})e_i - m^{-1}e_ie_je_i .\\
\end{eqnarray*}

We next show (\ref{eiejei-eq}).
Multiplying (R2) for $e_j$ by the left
and by the right with $e_i$, we find
$e_ie_jg_ie_je_i = le_ie_je_i$. Using (\ref{lm4}) we obtain
$(l+m^{-1})e_i - m^{-1}e_ie_je_i = le_ie_je_i$,
whence
$(l+m^{-1})e_ie_je_i = (l+m^{-1})e_i$.
As $lm\ne -1$, we find $e_ie_je_i = e_i$. This proves (\ref{B4}).

In order to prove the second equality of (\ref{lm1}),
we expand $g_ig_je_i$ by substituting the relation (\ref{B4}). We find
$$g_ig_je_i  =  g_ig_je_ie_je_i 
 =  e_jg_ig_je_je_i 
 =  l^{-1}e_jg_ie_je_i 
 =  e_je_i. 
$$

The first parts of the equalities of (\ref{7-eq}) and (\ref{8-eq}) are direct consequences of (\ref{lm1})
and (\ref{B4}).
In order to show the second part of
(\ref{8-eq}), we use the second equality of (\ref{lm1}) and (\ref{g-sq-eq}):
\begin{eqnarray*}
e_je_ig_j & = & (e_jg_ig_j)g_j 
 =  e_jg_i(ml^{-1}e_j-mg_j+1) \\
& = & me_j-me_jg_ig_j+e_jg_i 
 =  m(e_j-e_je_i) + e_jg_i. \\
\end{eqnarray*}
The second part of (\ref{7-eq}) follows from this by the anti-involution
of Remark \ref{opLm}(i).

For the first part of (\ref{gjgiej-inv-eq}), as the $g_i$ and $g_j$ are
invertible this is $g_ig_je_ig_jg_i=e_j$.  By (\ref{lm1}) the left side is
$e_je_ie_j$ which is $e_j$ by (\ref{B4}).

Finally we derive the second part of
(\ref{gjgiej-inv-eq}).
\begin{eqnarray*}
g_je_ig_j & = & g_je_ie_je_ig_j 
 =  (m(e_j-e_ie_j) + g_ie_j)e_ig_j \\
& = & me_je_ig_j - me_ie_je_ig_j + g_ie_je_ig_j \\
& = & m(m(e_j-e_je_i) + e_jg_i) - me_ig_j + g_i(m(e_j-e_je_i)+ e_jg_i))
\\
& = & m^2e_j-m^2e_je_i + m(e_jg_i - e_ig_j + g_ie_j)-mg_ie_je_i+ g_ie_jg_i
\\
& = & g_ie_jg_i + m^2e_j-m^2e_je_i + m(e_jg_i - e_ig_j +
g_ie_j)\\
&& -m(m(e_i-e_je_i) + g_je_i) \\
& = & g_ie_jg_i + m^2e_j-m^2e_i + m(e_jg_i - e_ig_j + g_ie_j - g_je_i).
\end{eqnarray*}
\end{proof}

\np
The above identities suffice for a full determination of the
BMW algebra associated with the braid group on 3 braids.

\begin{Cor}
The BMW algebra of type $\A_2$ has dimension 15 and is spanned by the monomials
\begin{eqnarray*}
&&1,\\
&&g_1, g_2, e_1, e_2, \\
&&g_1g_2, g_1e_2, g_2g_1, g_2e_1, e_1g_2, e_1e_2, e_2g_1, e_2e_1, \\
&&g_1g_2g_1, g_1e_2g_1.
\end{eqnarray*}
\end{Cor}

\begin{proof}
Let $B$ be the BMW algebra of type $\A_2$.
Of the sixteen possible words of length $2$ the eight consisting of two
elements with the same index can be reduced to words of length $1$. For, by (D1)
$g_i^2$ can be written as a linear combination of $g_i$, $e_i$ and $1$ and
by (\ref{e-sq-eq}) $e_i^2$ is a scalar multiple of $e_i$. Finally, by
relation (R1) the remaining four words reduce to $e_i$.

Now consider words of length $3$. By the knowledge that $x^{-1}e_i$ is
an idempotent and relation (\ref{B4}) it is clear that no words of length 3 can
occur containing only $e$'s. Words containing only $g$'s can be reduced if
two $g$'s with the same index occur next to each other. This leaves two
possible words $g_ig_jg_i$ either of which can be rewritten to the other one by (B1).

If a word contains $e$'s and $g$'s, no $e$ and $g$ may occur next to each
other having the same index as this can be reduced by relation (R1). So the
only sequences of indices allowed here are ${i,j,i}$ and ${j,i,j}$. If a $g$ occurs in the
middle, we can reduce the word by relation (R2) or (\ref{lm1}).
This leaves the case with an $e$ in the middle. By 
(\ref{8-eq}), (\ref{7-eq}), and (\ref{eiejei-eq}) these words reduce unless both the other
elements are $g$'s. Finally by (\ref{gjeigj-eq}) the two words
left, viz.~$g_ie_jg_i$ and $g_je_ig_j$, are equal up to some terms of shorter length, so at most one is in the
basis.

All words of length $4$ that can be
made by multiplication with a generator from the two words left of length $3$, can be reduced. First consider
$g_ig_jg_i$. Multiplication by a $g$ gives, immediately or after applying (B2), a
reducible $g^2$ component. Similarly, multiplication by an $e$ will result in a
reducible $e_ig_i$ word part.
This leaves us with multiples of $g_ie_jg_i$. As noted above, 
they can be expressed as a linear combination of $g_je_ig_j$ and terms of
shorter length.  Again, multiplication by $g$ leads to a $g^2$ component and
the word can be reduced. Multiplication by $e$ will always enable application
of relation (R2) to the constructed word and can therefore be reduced, proving
that no reduced words of length $4$ occur in $B$.

Finally, by use of the $15$ elements as a basis, one can construct an algebra
satisfying all relations of the BMW algebra, so the dimension of $B$ is
indeed $15$. This is done in \cite{Wenzl} and later in this paper.
\end{proof}

\begin{Prop} \label{lm0}\label{ij=ji-id-lm}
The following identities hold for $i\not\sim j$.
\begin{eqnarray}
e_ig_j &= & g_je_i,\\
\label{ijCommuteEq}
e_ie_j &= & e_je_i.
\end{eqnarray} 
\end{Prop}

\begin{proof}
By (D1), the $e_i$ are defined as polynomials in $g_i$ and belong to the
subalgebra of $B$ generated by $g_i$. By (B1) this subalgebra commutes
with $g_j$. 
\end{proof}

\begin{Prop}
There is a unique semilinear automorphism of $B$ of order 2 determined by
$$g_i\mapsto -g_i^{-1},\quad e_i\mapsto e_i, \quad l\mapsto -l^{-1},\quad m\mapsto m.$$
It commutes with the opposition involution
of Remark \ref{opLm}(i).
\end{Prop}

\nl
\begin{proof}
Using the identities proved above, it is readily verified that the
defining relations of $B$ are preserved.
\end{proof}

\np
We recall the definition of the BMW algebra as given in \cite{Wenzl}; however,
we take the parameters $q$, $r$ to be indeterminates over the field.  

\begin{Def} \label{wenzl}
  Let $q,r$ be indeterminates. The Birman-Murakami-Wenzl algebra $BMW_k$ is
  the algebra over $\C(r,q)$ generated by $1,g_1, g_2, \ldots, g_{k-1}$, which
  are assumed to be invertible, subject to the relations
\begin{eqnarray*}
g_ig_{i+1}g_i &=& g_{i+1}g_{i}g_{i+1},\\
g_ig_j &=& g_jg_i \mbox{ if } \mid i-j \mid \geq 2,\\
e_ig_i  &=& r^{-1}e_i \\
e_ig_{i-1}^{\pm 1}e_i &=& r^{\pm 1} e_i,\\
\end{eqnarray*}
where $e_i$ is defined by the equation
$ (q-q^{-1})(1-e_i) = g_i-g_i^{-1}$.
\end{Def}

\np
We now show that our definition of the BMW algebra of type $\A_n$
coincides with this one.

\begin{Thm}\label{wenzl-thm}
Let $n\ge 2$. The BMW algebra $B$ of type $\A_{n-1}$ is 
the Birman-Murakami-Wenzl algebra
$BMW_n$ where $l = r$ and $m = q^{-1} -q$.
\end{Thm}

\begin{proof}
To show both definitions are of the same algebra, we take our parameters $l
= r$ and $m = q^{-1}-q$. The first two relations for both algebras are
the same. It is
evident from the definition of $e_i$ in both $BMW_n$  and $B$
that $g_i$ and $e_i$ commute, so the third relation for $BMW_n$ is equivalent
to (\ref{eigi-eq}) and (R1) for $B$.  
Also the relation $e_ig_{i-1}e_i=le_i$ for $BMW_n$ is equivalent to (R2) for $B$.
To see that $g_i$ and $e_i$ in $B$ satisfy $e_ig_{i-1}^{-1}e_i=l^{-1}e_i$,
the final defining relation for $BMW_n$,
observe that, for $i\sim j$, 
by (\ref{g-inv-eq}), (R2), (\ref{e-sq-eq}), (\ref{B4}), and (\ref{def-x-eq}),
\begin{eqnarray*}
e_ig_j^{-1}e_i &=& e_i(g_j+m-me_j)e_i = (l+mx-m)e_i = l^{-1}e_i.
\end{eqnarray*}
The definition of $e_i$ follows from Remark \ref{opLm}(ii).  This shows that $B$
is a homomorphic image of $BMW_n$.  To
go the other way it is shown in \cite{Wenzl} (4) that $e_ig_{i+1}^{\pm
  1}e_i=r^{\pm 1}e_i$ and so all the relations of $B$ are
verified for $BMW_n$ except (D1).  This follows from
(10) in \cite{Wenzl} which when corrected reads
$g_i^2=(q-q^{-1})(g_i-r^{-1}e_i)+1$.  The invertibility of the $g_i$ follows
from (\ref{g-inv-eq}).  This shows the algebras are isomorphic.

\end{proof}

Although it is not needed for our computations, there is a cubic relation which is sometimes 
instructive.  
\begin{Prop}\label{cubicrel}  The elements $g_i$ of $B$ satisfy the cubic relation 
$$(g_i^2+mg_i-1)(g_i-l^{-1})=0.$$
\end{Prop}

\nl \begin{proof}
By (D1) and (\ref{eigi-eq}), we have
$(g_i^2+mg_i-1)(g_i-l^{-1}) = e_i(g_i-l^{-1})= 0$.
\end{proof}

In \cite{Wenzl} Proposition~$3.2$, it is shown that the algebras of type
$\A_{n-1}$, the so-called BMW algebras, are finite dimensional. This uses
in a crucial way that the symmetric group $\Sym_n\cong W(\A_{n-1})$ is doubly transitive on the cosets of
$\Sym_{n-1}$. This is not true for the other algebras. However, we
provide a proof of finite dimensionality which applies to the algebras
of type $\A_n$ as well.

Let $(W,R)$ be the Coxeter system of type $M$ and let $\{r_1,\ldots,r_n\}=R$.
Assume furthermore that $M$ is spherical.  Then the number of positive roots,
$|\Phi^+|$, is the length of the longest word in the generators $r_i$ of $W$.
This means that any product in $B$ of $g_i$ and $e_i$ of longer length can be
rewritten by using the relations (B1) and (B2) until one of $g_i^2$, $g_ie_i$,
$e_ig_i$, $e_i^2$ occurs as a subproduct for some $i$. In the Coxeter group,
$r_i$ has order $2$ so we can remove the square and obtain a word of shorter
length. In our algebra, we can rewrite the four words to obtain a linear
combination of words of shorter length. This leads to the following result.

\begin{Prop}\label{findim-prop}
  If the diagram $M$ is spherical, then any word in the generators of $B$ of
  length greater than $|\Phi^+|$ in $g_i$, $g_i^{-1}$, $e_i$ can be expressed
  as a sum of words of smaller length by using the defining relations of $B$.
  In particular, $B$ is finite dimensional.
\end{Prop}

\nl
\begin{proof} 
We can express $g_i^{-1} $ by $e_i$ and $g_i$ to get 
sums of words in $g_i$ and $e_i$. Suppose $w$ is a word in $g_i$ and
$e_i$ of length greater than 
$|\Phi^+|$. Consider the word in the Coxeter group $w'$ in $r_i$ where
each $g_i$, $e_i$ in $w$ is 
replaced by $r_i$. Notice that if $i\not \sim j$ that both 
$r_i$ and $r_j$ commute and that both $e_i$ and $g_i$ commute 
with both $e_j$ and $g_j$. In particular, the same changes can 
be made without changing $w$ or $w'$. Suppose the relation 
(B2) is used in $w'$,
$r_jr_ir_j=r_ir_jr_i$. Consider 
the same term in $w$ where $r_i$ are replaced 
by $g_i$, or $e_i$ and the same for $r_j$. 
We showed in the previous sections that all possible ways of replacing
the $r_i$ and $r_j$ by $e$ and $g$ elements reduces the word
except for $g_ig_jg_i=g_jg_ig_j$ and $g_ie_jg_i=g_je_ig_j + \omega$,
where $\omega$ is a linear combination of monomials of degree less than
3.
In fact they give words of length $2$ or, in the case $e_jg_i^{\pm 1}e_j$,
length $1$. 
If we arrive at $e_ig_i=g_ie_i$ we can replace it by (R1) with 
$l^{-1}e_i$ of shorter length. If we arrive at $g_i^2$ we use
(\ref{g-sq-eq})
to express it as a sum of words with $g_i^2$ replaced with 
$e_i$, $g_i$, and the identity. The same holds for $g_i^{-2}$ using 
the definition. If we arrive at $e_i^2 $ we 
can replace it with a multiple of itself. In all cases we 
can reduce the length. 

It is now clear that any word in $g_i$, $e_i$ can be written 
as a sum of the words of length at most $|\Phi^+|$ in $g_i$, $g_i^{-1}$, and 
$e_i$. 
\end{proof} 

\nl
{\bf Proof of Theorem \ref{findim-thm}}.
This is a direct consequence of the above proposition.

\section{Artin group properties}\label{sec:artin}
In this section, $M$ is a connected, simply laced, spherical Coxeter
diagram. This means 
$M = \A_n$ $(n\ge 1)$, $\D_n$ $(n\ge 4)$, or $\E_n$ $(n\in \{6,7,8\})$. We
shall often abbreviate this condition by writing $M\in\ADE$.

We let $(A,S)$ be an Artin system of type $M$, that is, a pair consisting of
an Artin group $A$ of type $M$ with distinguished generating set
$\{s_1,\ldots,s_n\}$ corresponding to the nodes of $M$. Similarly, we let
$(W,R)$ be the Coxeter system of type $M$, where $R$ is the set of fundamental
reflections $r_1,\ldots,r_n$.  We shall write $\Phi$ for the root system
associated with $(W,R)$ and $\Phi^+$ for the set of positive roots with
respect to simple roots $\a_1,\ldots,\a_n$ whose corresponding reflections are
$r_1,\ldots,r_n$.  There is a map $\psi: W\to A$ sending $x$ to the element
$\psi(x) = s_{i_1}\cdots s_{i_t}$ whenever $x = r_{i_1}\cdots r_{i_t}$ is an
expression for $x$ as a product of elements of $R$ of minimal length.  For
$\b\in\Phi$, we shall denote by $r_\b$ the reflection with root $\b$
and by $s_\b$ its image $\psi(r_\b)$ in $A$.  For a subset $X$ of $W$ we write
$\psi(X)$ to denote $\{\psi(w)\mid w\in X\}$.  The map $\psi$ is a section of
the morphism of groups $\pi: A\to W$ determined by $s_i\mapsto r_i$, that is,
$\pi\circ \psi$ is the identity on $W$.

Let $B$ be the BMW algebra of type $M$ over $\Q(l,x)$.  By means of the
composition of $\psi$ and the morphism of groups $A\to B^\times$, we find a
map $W\to B$. We shall write $\wh{w}$ or, if $r_{i_1}\cdots r_{i_t}$ is a
reduced expression for $w$, also $\wh{i_1\cdots i_t}$ to denote the image in
$B^\times$ of $w$ under this map.
In particular, $g_i = \wh{r_i} =\wh{i}$.

Let $g\in A$. By $g^{-op}$ we denote the anti-involution $op$ of $B$ introduced in
Remark \ref{opLm}(i) applied to the inverse of the image of $g$ in $B$, which is the same as
the inverse of the anti-involution applied to $g$, viewed as an element of $B$.

\begin{Lm}\label{ij-paths}
  Let $i$, $j$ be nodes of $M$. There is a unique element of minimal length in $W$, denoted by
  $w_{ji}$, such that $w_{ji} r_j w_{ji}^{-1} = r_i$.  It has the following
  properties.
\begin{enumerate}[(i)]
\item If $i=i_1\adj i_2\adj \cdots \adj i_q = j$ is the geodesic in $M$ from
  $i$ to $j$, then $\wh{w_{ji}} = \wh{i_{q-1}}\wh{i_q}\wh{i_{q-2}}\wh{i_{q-1}}\cdots\wh{ i_2}\wh{i_3}\wh{i_1}\wh{i_2}$.
\item $w_{ij}^{-1} = w_{ji}$.
\item $\wh{w_{ij}}^{op} = \wh{w_{ji}}$.
\item $\wh{w_{ij}}e_i = e_{j}e_{i_{q-1}}\cdots e_{i_2} e_{i} =
  e_j\wh{w_{ij}}$.
\item $\wh{w_{ij}}e_i = \wh{w_{ij}}^{-op}e_i = \wh{w_{ji}}^{-1}e_i$.
\end{enumerate}
\end{Lm}

\begin{proof}
  Consider the graph $\Gamma$ whose nodes are the elements of $\Phi^+$ and in
  which two nodes $\alpha$, $\beta$ are adjacent whenever there is a node $i$
  of $M$ such that $r_i\alpha = \beta$.  An expression $w = r_{i_1}\cdots
  r_{i_t}$ of an element $w$ of $W$ satisfying $w r_j {w}^{-1} = r_i$
  represents a path $\beta, r_{i_t}\beta, \ldots, r_{i_2}\cdots r_{i_t} \beta,
  w \beta = \alpha$ from $\beta$ to $\alpha$ in $\Gamma$. Clearly, if $w$ is
  of minimal length then this path is a geodesic. This geometric setting
  readily leads to a proof of (i).

\nl A geodesic in $\Gamma$ from $\alpha$ to $\beta$ is given by a backwards
traversal of the geodesic from $\beta$ to $\alpha$. The corresponding
element of $W$ is $w^{-1}$, whence (ii) and (iii).

\nl Finally, (iv) and (v) follow by induction from (i) and,
respectively, (\ref{gjgiej-eq})
and (\ref{7-eq}).
\end{proof}

\np
For a positive root $\b$, we write $\het(\b)$ to denote its height,
that is, the sum of its coefficients with respect to the $\a_i$.
Furthermore, 
$\Supp(\b)$ is the set of $k\in\{1,\ldots,n\}$ such that
the coefficient of $\a_k$ in $\b$ is nonzero.

\begin{Prop}\label{wbeta-prop}
For each node $i$ of $M$ and each positive root $\beta$ 
there is a unique element $w\in W$ of minimal length such that $w\alpha_i =
\beta$. This element satisfies the following properties:
\begin{enumerate}[(i)]
\item If $\beta=\alpha_j$ for some $j$, then $w= w_{ij}$.
\item If $j$ is the unique node of $M$ in
$\Supp(\beta)$ nearest to $i$, then $l(w) = \het(\beta) + l(w_{ij}) -1$.
\end{enumerate}
\end{Prop}

\begin{proof}
Suppose first that $i$ lies in the support of $\beta$. Then $\beta$ can be
obtained from $\alpha_i$ by building up with addition of one fundamental
root at a time, which corresponds to finding an element $w$ of $W$ by
multiplication to the right of the fundamental reflection corresponding to
the newly added fundamental root. This shows that there exists $w\in W$ of
length at most $\het(\beta) - 1$ such that $w\alpha_i = \beta$. But the
height of $\beta$ is clearly at most $l(w)+1$, so the minimal length of any
element $w$ of $W$ so that $w\alpha_i = \beta$ must be $\het(\beta)-1$.

Next suppose that $i$ does not lie in the support of $\beta$ and let $j$ be
the nearest node to $i$ in the support of $\beta$. Then, with $y\in W$ as in
the first paragraph with respect to $\beta$ and $j$ so that
$y\alpha_j = \beta$ and $l(y) = \het(\beta)-1$, we have that $y w_{ji}\alpha_i
= \beta$
and that $l(yw_{ij})\le l(w) + l(w_{ij}) = \het(\beta)+l(w_{ij}) -1$.
On the other hand, in order to transform $\alpha_i$ into $\beta$ by a chain of
roots
differing by a fundamental root, we need to apply each root but $i$ and $j$ on the
geodesic in $M$ from $i$ to $j$ at least twice (once for creation of the
presence of the node in the support, and one for making it vanish). We also
need both $i$ and $j$ at least once. Hence, in order to make a fundamental
root
of the support of $\beta$ occur in the image $u\alpha_i$ of $\alpha_i$ of some
$u\in W$, we need
$l(u)\ge l(w_{ij})$, with equality only if $u = w_{ij}$ and $u\alpha_i =
\alpha_j$.
Notice that the fundamental reflections in $w_{ij}$ except for $\a_j$
do not contribute at all to the creation of the fundamental nodes in the
support of $\beta$, so that the estimate for the fundamental roots needed to
build
up $\beta$ stays as before. 
Taking $w = yu$ we find
$l(w) = l(yu) = l(y) + l(u) = l(y)+l(w_{ij}) = \het(\beta)+l(w_{ij}) -1$.

Next we prove uniqueness of $w$ as stated. Suppose $v\in W$ also satisfies
$l(v) = \het(\beta)+l(w_{ij}) -1$. As argued above, we must have
$v = v' w_{ji} $ and $l(v) = l(v') + l(w_{ji})$ so, without loss of
generality,
we may assume $i=j$ lies in the support of $\beta$.
If $l(w) = 0$ then there is nothing to show.
Suppose therefore $l(w)>0$ and apply induction on $l(w)$.
Take nodes $k,h$ of $M$ such that $l(r_kw)<l(w)$ and $l(r_hv)<l(v)$ while
$r_k\beta = \beta - \alpha_k$ and $r_h\beta = \beta - \alpha_h$.
Such $k$ and $h$ exist by the way $\beta$ is built up of fundamental roots
via $w$ and $v$, respectively.
Notice that 
$(\beta,\alpha_k) = (\beta,\alpha_h) = 1$.
Now consider $(\beta-\alpha_k,\alpha_h)$.
The value equals $-1$ if $k=h$;
$1$ if $h\ne k\not\sim h$; and $2$ if $k\adj h$.
In the first case, we apply induction to $(r_hw)\alpha_i = \beta-\alpha_h =
(r_hv)\alpha_i$,
and find $r_hw = r_hv$, whence $w = v$.

In the non-adjacent case, $\beta-\alpha_h-\alpha_k$ is also a root, so there
is a unique minimal $u\in W$ such that $u\alpha_i = \beta-\alpha_h-\alpha_k$.
Now $r_hr_ku\alpha_i = \beta = w\alpha_i = v\alpha_i$, so $r_hw\alpha_i = r_ku\alpha_i$
and $r_kv\alpha_i = r_hu\alpha_i$, whence, by
induction, both $r_hw = r_ku$ and $r_hu = r_kv$. But then
$w = r_hr_ku =r_kr_hu = v$.

Finally, if $k\adj h$, we find 
$(\beta-\alpha_k,\alpha_h)=2$, whence
$\beta = \alpha_h+\alpha_k$. But then $i$ must be either $h$ or $k$.
Assuming (without loss of generality) $i=h$, we find
$w=r_k$ and $v = r_h=r_i$, a contradiction with $v\alpha_i = \alpha_i+\alpha_h$. 

This establishes that $w$ is unique, and finishes the proof of the lemma.
\end{proof}

\begin{defn}
For a node $i$ of $M$ and a positive root $\beta$ 
we denote by $w_{\beta,i}$ the unique element (by the above proposition) of minimal length in $W$ for which
$w_{\beta,i}\alpha_i = \beta$.
We denote by $D_i$ the set $\{w_{\beta,i}\mid \beta\in\Phi^+\}$.
\end{defn}

If $w\in D_i$ then $w r_i w^{-1}$ is a shortest expression of
the reflection corresponding to $w\alpha_i$ as a conjugate of $r_i$.

\begin{Cor}\label{Li-props-lem}
For each node $i$ of $M$,
the  set $D_i$ satisfies the following properties, where $j$ is a node of $M$.
\begin{enumerate}[(i)]
\item If $r_jv\in D_i$ and $v\in W$ with $l(r_jv) = l(v)+1$,
then $v\in D_i$.
\item $w_{ij}\in D_i$.
\end{enumerate}
\end{Cor}

\begin{Lm}
If $i$ and $j$ are nodes of $M$, then $\wh{w_{\a_j,i}}e_i = \wh{w_{ij}}e_i$.
\end{Lm}

\begin{proof}
Building up 
$w_{\a_j,i}$ from the right, and letting the intermediate results act on
$\a_i$,
we find a shortest path $i=i_1\adj i_2\adj\cdots \adj i_t = j$ in $M$ from 
$i$ to $j$. The
element
$\wh{w_{ij}}$ represents the corresponding element $ \wh{i_{t-1}i_t} \cdots
\wh{i_2i_3}\wh{ii_2}$ of $B$.
\end{proof}

\begin{Lm}\label{conj-induction}
For all nodes $i,j,k$ of $M$ we have
$\wh{w_{ki}}\wh{w_{jk}}e_j = \wh{w_{ji}}e_j$.
\end{Lm}

\begin{proof} 
Denote by $i=i_1\adj i_2\cdots\adj i_q =k$
the geodesic from $i$ to $k$ and by 
$k=k_1\adj k_2\cdots\adj k_p = j$
the geodesic from $k$ to $j$.
Then there is an $m\in\{1,\ldots,q\}$ such that
$k=k_1=i_q\adj k_2 = i_{q-1}\adj\cdots k_m = i_{q-m+1}$
and $k_{m+1}\ne i_{q-m}$.
Then the geodesic from $i$ to $j$ is
$i=i_1\adj i_2\cdots\adj i_{q-m}\adj k_m\adj k_{m+1}\adj\cdots\adj
{k_{p-1}}\adj k_p$
and so
\begin{eqnarray*}
\wh{w_{ki}}\wh{w_{jk}}e_j &=& 
\wh{w_{ki}} e_{k_1}\cdots e_{k_p}\\
&=&
e_{i_1}\cdots e_{i_q} e_{k_1}\cdots e_{k_p}\\
&=& e_{i_1}\cdots e_{i_{q-m}}e_{k_m} \cdots 
e_{k-1} e_ke_{k-1}\cdots e_{k_p}\\
&=& e_{i_1}\cdots e_{i_{q-m}}e_{k_m} \cdots 
e_{k-1} e_k e_{k-1}\cdots e_{k_{m}}e_{k_{m+1}}\cdots e_{k_p}\\
&=& e_{i_1}\cdots e_{i_{q-m}}e_{k_m} e_{k_{m+1}}\cdots  e_{k_p}\\
&=& \wh{w_{ji}} e_j .
\end{eqnarray*}
\end{proof}

\np
For $\a,\b\in\Phi^+$ with $\a \le \b$ (that is, for each $i$, the difference
of the coefficient of $\a_i$ in $\b$ and the coefficient of $\a_i$ in $\a$ is
nonnegative), let $w_{\b,\a}$ be the (unique) shortest element of $W$ mapping
$\a$ to $\b$. Clearly, $l(w_{\b,\a}) = \het(\b)-\het(\a)$.  Thus, $w_{\b,i} =
w_{\b,\a_i}$ if $i\in\Supp(\b)$.  For a positive root $\b$, set
$d_\b=\psi(w_{\a_0,\b}^{-1})\in A$. This implies that $s_{\a_0} = d_\b^{op}
s_{\b} d_{\b}$.  For a node $i$ such that $\a_i$ is orthogonal to $\b$, we
shall need the following Artin group element.
\begin{eqnarray}\label{hbetai-df}
h_{\b,i} &=& d_{\b} ^{-1}s_i d_{\b}.
\end{eqnarray}

\begin{Lm}\label{H-lem}
The following relations hold for elements $h_{\c,k}$ of the Artin group $A$,
where we are always assuming that $\c$ is a positive root and $(\a_k,\c) = 0$.
\begin{eqnarray}
h_{\b,i}h_{\b,j} & = & h_{\b,j}h_{\b,i} \ \qquad \mbox{ if } 
i\not\sim j \label{H-com-rel}\\
h_{\b,i}h_{\b,j}h_{\b,i} & = & h_{\b,j}h_{\b,i}h_{\b,j} \
\qquad \mbox{ if } i\sim j \label{H-braid-rel}\\
h_{\b+\a_j,i} & = & h_{\b,i} \ \ \ \ \qquad \mbox{ if } i\not\sim j \label{H-orth-rel}\\
h_{\b+\a_j,i} & = & h_{\b-\a_i,j} \ \ \qquad \mbox{ if } 
i\sim j \label{H-pm-rel}\\
h_{\b-\a_i-\a_j,i} & = & h_{\b,j} \ \ \ \qquad \hfill \mbox{ if } 
i\sim j \label{H-10-rel} \\
h_{\b+\a_i+\a_j,j} & = & h_{\b,i} \ \ \ \qquad \hfill \mbox{ if } 
i\sim j\label{H-0-1-rel}\\
h_{\a_i,j} & = & h_{\a_j,i} \ \ \ \qquad \hfill \mbox{ if } 
i \mbox{ and } j \mbox{ are at distance 2 in } M \label{H-dist2-rel}\\
 h_{\a_j,k} &=& h_{\a_i,k} \ \ \ \qquad \hfill \mbox{ if } 
i \sim j. \label{H-jjprime-rel}
\end{eqnarray}
\end{Lm}

\begin{proof}
The rules are all straightforward applications of corresponding rules for
$d_\b$. We prove (\ref{H-pm-rel}) and (\ref{H-jjprime-rel}) and leave the rest 
to the reader.  

For rule (\ref{H-pm-rel}), we have $d_{\beta-\a_i} =
s_is_jd_{\b+\a_j}$ in the Artin group whereas $i\sim j$, $(\a_i,\b)=-1$, and
$(\a_j,\b) = 1$, so $h_{\b-\a_i,j}$ is the Hecke algebra element corresponding
to the Artin group element $d_{\b+\a_i}^{-1}s_jd_{\b+\a_i} =
d_{\b-\a_j}^{-1}s_j^{-1}s_i^{-1}s_js_is_jd_{\b-\a_j} =
d_{\b-\a_j}^{-1}s_id_{\b-\a_j}$, and so
$h_{\b-\a_i,j}$ coincides with $h_{\b-\a_j,i}$.

We finish with (\ref{H-jjprime-rel}). It is a direct consequence of $s_i^{-1}d_{\a_j}=
d_{\a_i+\a_j} =s_j^{-1}d_{\a_i}$
and the fact that $k$ is adjacent to neither $i$ nor $j$:
$$
h_{\a_j,k}= d_{\a_j}^{-1}  k d_{\a_j} = 
d_{\a_i}^{-1}s_js_{i}^{-1}  s_{k}s_is_j^{-1} d_{\a_j} = 
d_{\a_i}^{-1}  s_{k} d_{\a_i} =
 h_{\a_i,k} .
$$
\end{proof}

As before, let $C$ be the set of nodes $i$ of $M$ for which $\a_i$ is orthogonal to the highest root
$\a_0$ of $\Phi^+$.

\begin{Lm}\label{Hbeta-lm}
The following properties hold for $C$.
\begin{enumerate}[(i)]
\item If $i$ is a node of $M$ and $\b\in\Phi^+$ satisfies $(\a_i,\b) = 0$, then 
there is a node $j$ of $C$ such that $h_{\b,i} = s_j$.
\item For each $j$ in $C$ there exist non-adjacent nodes $i$, $k$ with $h_{\a_i,k} = s_j$.
\end{enumerate}
\end{Lm}

\begin{proof}
(i). If $\b=\a_0$, then $i$ is a node orthogonal to $\a_0$ and so $h_{\b,i} = s_i$
and $i$ belongs to $C$ by
definition of $C$.
We continue by induction with respect to the height of $\beta$.
Assume $\het(\b)<\het(\a_0)$. Then there is a node $j$ such that
$(\a_j,\b) = -1$, so $\c = \b+\a_j$ is a root, whence $d_\b = s_jd_\c$.
If  $i\not\sim j$, then, by (\ref{H-orth-rel}), $h_{\b,i}= h_{\c,i}$.
Otherwise, by (\ref{H-0-1-rel})
$h_{\b,i}= h_{\c+\a_i,j}$. In both cases the expression found for $h_{\b,i}$ is as
required by the induction hypothesis.

\nl(ii).
Let $j$ be a node in $C$. Then $h_{\a_0,j} = \wh{j}$.
Let $\b$ be a minimal positive root for which there exists
a node $k$ with $(\a_k,\b) = 0$ and $h_{\b,k} = \wh{j}$.
If $\het(\b) > 1$, 
take a node $i$ such that $(\a_i,\b) = 1$.
By Lemma \ref{H-lem}, either $i\sim k$ and
$h_{\b-\a_i-\a_k,i} = \wh{j}$, or $(\a_i,\a_k) = 0$ and $h_{\b-\a_i,k} = \wh{j}$.
Therefore, we may assume $\het(\b) = 1$, and so $\b=\a_i$ for some $i$
with $(\a_i,\a_k) = 0$.
\end{proof}

\begin{Lm}\label{isb=sbi-lm}
If $i$ is a node of $M$ and $\b$ a positive root such that $(\a_i,\b) = 0$, then
$$s_is_\b = s_\b s_i.$$
\end{Lm}

\begin{proof}
We proceed by induction on $\het(\b)$.
If $\het(\b) = 1$, then $\b=\a_j$.
As $(\a_i,\b) = 0$, we have $i\not\sim j$ and so $s_is_\b = s_is_j = s_js_i = s_\b s_i$ by the
braid relations.

Assume now that $\het(\b) >1 $.  Let $j$ be a node of $M$ such that $(\a_j,\b)=1$, 
so $\b-\a_j$ is a positive root.  Then $s_\b = s_js_{\b-\a_j}s_j$. If $j\not\sim i$, then
$(\a_i,\b-\a_j) = 0$, so, by the induction hypothesis,
$s_is_{\b-\a_j} = s_{\b-\a_j}s_i$, whence
$s_is_\b = s_is_js_{\b-\a_j}s_j = s_js_is_{\b-\a_j}s_j = s_js_{\b-\a_j} s_is_j
= s_js_{\b-\a_j} s_js_i = s_\b s_i$.
Otherwise, $j\sim i$, and $\c = \b-\a_i-\a_j$ is a positive root with
$(\a_j,\c) = 0$ and $s_\b = s_js_is_\c s_is_j$. By the induction hypothesis,
$s_js_{\c} = s_{\c}s_j$, whence
$s_is_\b = s_is_js_is_{\c}s_is_j = s_js_is_js_{\c}s_is_j = 
s_js_is_{\c}s_js_is_j = s_js_is_{\c}s_is_js_i = s_{\b}s_i$. 
\end{proof}

\section{Some ideals of the BMW algebra}
\label{sec:ideals}
In this section, let $M$ be a simply laced Coxeter diagram (not necessarily
spherical).  In the BMW algebra $B$ of type $M$, the $e_i$ generate an ideal (by
which we mean a 2-sided ideal).  Taking products of $e_i$'s for non-adjacent
nodes $i$ of $M$, we obtain ideals all contained in the ideal generated by $e_1$.   

\begin{defn}\label{I_jdefs}
Let $Y$ be a coclique of $M$, that is, a subset of the nodes of $M$
in which no two nodes are adjacent. The {\em ideal of type $Y$}
is the (2-sided) ideal of $B$ generated by $e_Y$, where
$$e_Y = \prod_{y\in Y} e_y.$$
The element $e_Y$ is well defined as the product does not depend on the order of
the $e_y$ in view of (\ref{ijCommuteEq}).
The ideal $Be_YB$ is denoted by $I_Y$.
By $I_j$, for $j=1,\ldots,n$, we denote the ideal generated by all $I_Y$ for
$Y$ a coclique of size $j$.
\end{defn}

\np
Since the $e_i$ are scalar multiples of idempotents, so are their products $e_Y$
for $Y$ a coclique of $M$.

\begin{Prop}\label{Ideal-I-prop}
Let $X$, $Y$ be cocliques of $M$. 
\begin{enumerate}[(i)]
\item If $X\subseteq Y$ then $I_Y\subseteq I_X$.
\item If $\{r_j\mid j\in X\}$ is in the same $W$-orbit as $\{r_j\mid j\in Y\}$ 
then $I_X = I_Y$.
\item
The quotient algebra $B/I_1$ is the Hecke algebra of type $M$ over
$\Q(l,x)$, with parameter $m$.
\end{enumerate}
\end{Prop}

\begin{proof}
(i) is immediate from the definition of $I_Y$ and the commutation
of the $e_i$ for $i\in Y$.

\nl (ii). For $|X| = |Y| = 1$, say $X=\{i\}$ and $Y = \{j\}$,
this follows from the existence of
the invertible element $\wh{w_{ij}}$ as in Lemma \ref{ij-paths}(iv).
More generally, by \cite{Godelle}, there exists $w\in W$ such that
$\wh{w}\wh{X}\wh{w}^{-1} = \wh{Y}$. This implies
$\wh{w}e_X\wh{w}^{-1} = e_{Y}$, whence $I_X = I_Y$.

\nl (iii).
By (\ref{lm1}),
invertibility of the $g_i$ and connectedness of $M$, the ideal $I_1$
coincides with $I_{\{j\}}$
for any node $j$ of $M$.  Consequently, the quotient ring $B/I_1$ is obtained by
setting $e_i = 0$ for all $i$.
This means that the braid relations (B1) and (B2) and
(D1) are the defining relations for $B/I_1$ in terms of $g_i$.
Now (D1) reads $g_i^2+mg_i-1=0$, so we obtain the defining relations
of the Hecke algebra.
\end{proof}

\np
By (i), we have the chain of ideals
$$I_1 \supset I_2 \supset \ldots \supset I_k,$$
where $k$ is the maximal
coclique size of $M$.  By analogy with the BMW algebra of type $\A_n$ and
computer results for $\D_4$ we expect this is a strictly decreasing series of
ideals.  We already know from (iii) of the above proposition that $I_1$ is
properly contained in $B$. Straightforward calculations for the
Lawrence-Krammer representation, described in \cite{CW} and in \cite{paris}
for the non-spherical types, show that (D1), (R1), (R2) are also satisfied and
so that it is a representation of $B$. Furthermore it can be seen that $e_i$
is not represented as $0$ but $e_ie_j$ is for any two distinct non-adjacent
nodes $i$, $j$ of $M$.  These calculations will be presented in a more general
setting later, in Section \ref{sec:krammer}.  As a consequence $I_2$ is
properly contained in $I_1$. This follows also of course from Theorem
\ref{main-thm}.

\np
It is also clear from the definition that $I_j=\{0\}$ when $j$ is bigger than
the maximal coclique size of $M$.
These numbers are
$\lfloor (n+1)/2\rfloor$ for $\A_n$;
$\lfloor n/2\rfloor +1 $ for $\D_n$;
$3$ for $\E_6$; and 4 for both $\E_7$ and $\E_8$.

\section{Structure of $I_1/I_2$}\label{sec:I1modI2}
Throughout this section, $M$ is a connected simply laced spherical diagram.
This means $M\in \ADE$. By $B$ we denote the corresponding BMW algebra over
$\Q(l,x)$, by $(A,S)$ the corresponding Artin system, and by $(W,R)$ the
corresponding Coxeter system.  Furthermore, $\Phi^+$ is the set of positive
roots associated with $(W,R)$ and $C$ the set of nodes $i$ of $M$ with $\a_i$
orthogonal to the highest root of $\Phi^+$.

We now prepare for considerations of $B$ modulo $I_2$. Some of the results
hold for $B$ and others just modulo $I_2$. This is indicated in the statements.
The aim is to find a linear spanning set for $I_1/I_2$ of size
$|\Phi^+|^2|W_C|$.  In particular, we
obtain an upper bound for $\dim(I_1/I_2)$, which by Theorem \ref{main-thm}
will be an equality.

\np Let $i$ be a node of $M$ and let $Z_i$ be the subalgebra (not necessarily
containing the identity) of $B$ generated by all elements of the form
$\wh{w_{ji}}\wh{k}\wh{w_{ij}}e_i$ for $j$ and $k$ non-adjacent nodes of $M$.
We allow for $j$ and $k$ to be equal, so that, in case $M = \A_2$,
the subalgebras $Z_i$ are one-dimensional (scalar multiples of $e_i$).
By Lemma \ref{ij-paths}(iv),(v), the generators can be written in various ways:
$$e_i\wh{w_{ji}}\wh{k}\wh{w_{ji}}^{-1} =
\wh{w_{ji}}\wh{k}\wh{w_{ji}}^{-1}e_i
=\wh{w_{ji}}\wh{k}\wh{w_{ij}}e_i.$$

We will need an integral version of $Z_i$ and $B$.
We shall work with the coefficient ring $E = \Q(x)[l^{\pm}]$ inside
our field $\Q(l,x)$. Observe $m\in E$ by (\ref{def-x-eq}).
Let $B^{(0)}$ be the subalgebra of $B$ over $E$ generated by all $g_i$ and
$e_i$, and let $Z_i^{(0)}$ be the subalgebra of $Z_i$ over
$E$ generated by the same elements as taken above for generating $Z_i$.
Then $Z_i^{(0)}$ is a subalgebra of $B^{(0)}$.

\begin{Prop}\label{Z_i-prop}
The subalgebra $Z_i^{(0)}$ of $B^{(0)}$ satisfies the following properties.
\begin{enumerate}[(i)]
\item It centralizes $e_i$ and has identity element $x^{-1}e_i$.
\item $Z_i^{(0)} = \wh{w_{ji}}Z_j^{(0)}\wh{w_{ji}}^{-1}$ for all nodes $j$ of $M$.
\item The scaled versions
$x^{-1}e_i\wh{w_{ji}}\wh{k}\wh{w_{ji}}^{-1}$  of the generators of $Z_i^{(0)}$
satisfy the quadratic relation
$X^2+mX-1_i = 0 \mod I_2$, where $1_i$ stands for the identity element
$x^{-1}e_i$ of $Z_i^{(0)}$.
\end{enumerate}
\end{Prop}

\begin{proof}
(i). Since $x^{-1}e_i$ is an idempotent (cf.~(\ref{e-sq-eq})),
it suffices to verify that the generators of $Z_i$ centralize $e_i$. This
follows
from the following computation, in which Lemmas
\ref{ij-paths} and \ref{conj-induction} are used.
$$
\wh{w_{ji}}\wh{k}\wh{w_{ji}}^{-1} e_i
=
\wh{w_{ji}}\wh{k}e_j \wh{w_{ij}} 
=
\wh{w_{ji}}\wh{k}e_j \wh{w_{ji}}^{-1} 
=
\wh{w_{ji}}e_j \wh{k} \wh{w_{ji}}^{-1} 
= e_i\wh{w_{ji}} \wh{k} \wh{w_{ji}}^{-1}.
$$

\nl(ii). 
For the generator $e_h\wh{w_{jh}}\wh{k} \wh{w_{jh}}^{-1}$ of $Z_h^{(0)}$,
where $j\perp k$, we have
\begin{eqnarray*}
\wh{w_{hi}}e_h\wh{w_{jh}}\wh{k} \wh{w_{jh}}^{-1}\wh{w_{hi}}^{-1} 
&=&
\wh{w_{hi}}\wh{w_{jh}}e_j\wh{k} \wh{w_{jh}}^{-1}\wh{w_{hi}}^{-1} 
=
\wh{w_{ji}}e_j\wh{k} \wh{w_{jh}}^{-1}\wh{w_{hi}}^{-1} \\
&=&
\wh{w_{ji}}\wh{k}e_j \wh{w_{jh}}^{-1}\wh{w_{hi}}^{-1} 
=
\wh{w_{ji}}\wh{k}e_j \wh{w_{hj}}\wh{w_{hi}}^{-1} \\
&=&
\wh{w_{ji}}\wh{k} \wh{w_{hj}}e_h\wh{w_{hi}}^{-1} 
=
\wh{w_{ji}}\wh{k} \wh{w_{hj}}e_h\wh{w_{ih}} \\
&=&
\wh{w_{ji}}\wh{k} e_j \wh{w_{hj}}\wh{w_{ih}} 
=
\wh{w_{ji}}\wh{k} e_j \wh{w_{ij}} \\
&=&
\wh{w_{ji}}\wh{k}e_j \wh{w_{ji}}^{-1} 
=
\wh{w_{ji}}e_j \wh{k}\wh{w_{ji}}^{-1} \\
&=&
e_i\wh{w_{ji}}\wh{k}\wh{w_{ji}}^{-1} ,
\end{eqnarray*}
whence 
$\wh{w_{hi}}Z_h^{(0)}\wh{w_{hi}}^{-1} \subseteq Z_i^{(0)}$. The rest follows easily.

\nl (iii).
Substituting $x^{-1}e_i\wh{w_{ji}}\wh{k}\wh{w_{ji}}^{-1}$ for $X$, we find
\begin{eqnarray*}
&&(x^{-1}e_i\wh{w_{ji}}\wh{k}\wh{w_{ji}}^{-1} )^2
+m (x^{-1} e_i\wh{w_{ji}}\wh{k}\wh{w_{ji}}^{-1}) -x^{-1}e_i \\
&&\qquad = x^{-1}e_i\wh{w_{ji}}(\wh{k}^2+m\wh{k}-1) \wh{w_{ji}}^{-1}   =
x^{-1}e_i\wh{w_{ji}}e_k \wh{w_{ji}}^{-1} \in Be_je_kB\subseteq I_2.
\end{eqnarray*}
\end{proof}

\np
We recall that $w_{\b,i} \in W$ is
the element of minimal length with the property that $w_{\b,i} \a_i = \b$
with
$\a_i,\b \in \Phi^+$.

\begin{Lm}\label{ai+aj}  Suppose $i$, $j$, and $k$ are distinct nodes of $M$.  Then 
$$
e_i\wh{j}e_k = 
\left\{
\begin{tabular}{rcl}
$e_ie_k\wh{j}$&if& $j\not \sim k$ and $i \not \sim k$,\\
$\wh{w_{\a_i,k}}e_k\wh{j} $ &if& $j\not \sim k$ and  $i \sim k$,\\
$\wh{w_{\a_i,k}}e_k(\wh{i} + m)-me_ie_k$  &if& $j\sim k$ and $i\sim j$,\\
$\wh{w_{\a_i,k}}e_k \wh{jkikj}$&if&$j\sim k$ and $i \sim k$,\\
$e_ie_k\wh{w_{ik}}\wh{j}\wh{w_{ki}}$  &if& $j\sim k$, $i\not \sim j$, and $i\not \sim k$.
\end{tabular}
\right.
$$
In each case the result is in $\wh{w_{\a_i,k}}Z_k^{(0)} + I_2$.
\end{Lm}

\begin{proof}
In the first two cases as $j\not \sim k$ we have $e_i\wh{j}e_k=e_ie_k\wh{j} $.  If 
$i\not \sim k$, $e_ie_k$ is in $I_2$.  If $i\sim k$, $e_ie_k=w_{\a_i,k}e_k$.  These 
are the only possibilities when $j\not \sim k$.

Suppose next that $j\sim k$.  In the last case $e_i$ commutes with $\wh{j}$ and 
$e_ie_k$ is in $I_2$.  Suppose then $i \adj j$. Of course then $i \not\adj k$ since the type is spherical.
Now
by (R2)
\begin{eqnarray*}
e_i \wh{j} e_k & = & (e_i e_j e_i) \wh{j} e_k 
 =  (e_i e_j \wh{ij}) \wh{j} e_k \\
& = & e_i e_j \wh{i} (1 - m\wh{j} +ml^{-1}e_j) e_k 
 =  e_i e_j \wh{i} e_k - m e_i e_j e_i e_k + m e_i e_j e_k \\
& = & e_i e_j e_k (\wh{i} + m) - m e_i e_k 
= \wh{w_{\a_i,k}}e_k  (\wh{i} + m) - me_ie_k
\end{eqnarray*}
As $e_ie_k\in I_2$ the result follows. 

\np
Finally, if $i \adj k$ then necessarily $i \not\adj j$, and
$$
e_i \wh{j} e_k  =  
e_ie_ke_i \wh{j} e_k = 
e_ie_k \wh{j}e_i e_k = 
e_ie_k e_i \wh{j} \wh{ki} = 
e_ie_k \wh{ikjki} =
e_i e_k \wh{jkikj}
$$
In each of the cases the elements are in $\wh{w_{\a_i,k}}Z_k^{(0)}+I_2$ from the 
definition.   
\end{proof}

\np
If some of $i$, $j$, $k$ are equal, similar results follow from the defining
relations
and Propositions \ref{iji=jij-id-lm} and \ref{lm0}.

\begin{Lm} \label{H0}
Let $i,j,k\in \{1,\ldots,n\}$
and let $\gamma$ be the shortest path from $j$ to $k$.
Then 
$$
\wh{i} \wh{w_{\a_j,k}} e_k = 
\left\{
\begin{tabular}{rcl}
$\wh{w_{\a_j,k}}e_k \wh{i}$ &if& $i\not \sim $ any point of $\gamma$,\\
$l^{-1}\wh{w_{\a_j,k}}e_k$  &if& $i=j$, \\
$\wh{w_{\a_j,k}}e_k\wh{h'}\mod I_2$ &if& $i\in \gamma$, $i\ne j$, $h'$\\
                && on the path from $i$ to $j$,\\
                && $h'$ at distance 2 to $i$ in $M$,\\
$\wh {w_{\a_j,k}}e_k\wh{w_{h'k}}\wh{i}\wh{w_{kh'}} $ &if& $i\not \in
\gamma$, $i\sim h $, $h\in \gamma$, \\
&& $h\ne j$, $h'\sim h$, and \\
&&$h'$ on the path from $h$ to $j$,\\
$\wh{w_{\a_i+\a_j,k}}e_k + m\wh{w_{\a_i,k}}e_k-m\wh{w_{\a_j,k}}e_k $ &if& $i\in \gamma$ and  $i\sim j$,\\
$\wh{w_{\a_i+\a_j,k}}e_k  $ &if& $i\not \in \gamma$ and $i\sim j$.\\
\end{tabular}
\right.
$$
Also 
$$
e_i\wh{w_{\a_j,k}} e_k = 
\left\{
\begin{tabular}{rcl}
$x\wh{w_{\a_j,k}}e_k$  &if& i=j, \\
$0\mod I_2$    &if&  $i\not \sim j$, \\
$\wh{w_{\a_i,k}}e_k $ &if& $i\sim j$.  \\
\end{tabular}
\right.
$$
In each case, the result is in $\wh{w_{r_i\a_j,k}}Z_k^{(0)} +
m\wh{w_{\a_j,k}}Z_k^{(0)}
+ m\wh{w_{\a_i,k}}Z_k^{(0)} + I_2$.
\end{Lm}

\begin{proof}
Consider the shortest path $\gamma = k,\ldots,j$ from $k$ to $j$ in $M$.
If $i$ is non-adjacent to each element of this path, then
the statement holds.
Also if $i = j$ the statement follows immediately.
This leaves two possibilities, $i$ is in $\gamma$, or $i$ is not 
in $\gamma$ but is adjacent to some $h$ in $\gamma$.

Assume that $i$ occurs in $\gamma$.
If $i\sim j$, then
by (\ref{7-eq})
\begin{eqnarray*}
\wh{i} \wh{w_{\a_j,k}} e_k &=& \wh{i}e_je_i\cdots e_k\\
&=& \wh{j}^{-1}\wh{w_{ki}} e_k =
\wh{w_{\a_i+\a_j,k}} e_k
+m \wh{w_{ki}} e_k - m \wh{w_{kj}} e_k.
\end{eqnarray*}

Suppose, therefore, that $i\not\sim j$. Then $\wh{i} \wh{w_{kj}} e_k =
e_j \cdots 
e_{h'} \wh{i} e_h e_i e_{i'}\cdots e_k$ with $h' \adj h \adj i \adj i'$.
Substitution 
of $\wh{i} e_h e_i = \wh{h} e_i - m e_h e_i + m e_i$ and use of Lemma~\ref{H0} gives
\begin{eqnarray*}
\wh{i} \wh{w_{kj}} e_k &=& e_j \cdots e_{h'} \wh{i} e_h e_i \cdots e_k
= e_j \cdots e_{h'} (\wh{h} e_i - m e_h e_i + m e_i) e_{i'} \cdots e_k\\
&=& e_j \cdots e_{h'} \wh{h} e_i e_{i'} \cdots e_k - m \wh{w_{kj}} e_k +
m e_j 
\cdots e_{h'} e_i \cdots e_k \\
&\in & e_j \cdots e_{h'} e_h e_i (\wh{h'} +m )e_{i'} \cdots e_k -
m\wh{ w_{kj} }
e_k + I_2 \\
&=& e_j \cdots e_{h'} e_h e_i \wh{h'} e_{i'} \cdots e_k + I_2 
= e_j \cdots e_{h'} e_h e_{i} \cdots e_k\wh{h'} + I_2 \\
&=& \wh{w_{kj}} e_k h' + I_2 \\
\end{eqnarray*}

\np
Next assume $i$ is not in $\gamma$ but is adjacent to some $h$ in
$\gamma$. Suppose there exists $h'\adj h $ in $\gamma$, so
$$\wh{i} \wh{w_{kj}} e_k = e_j \cdots e_{h'} \wh{i} e_h \cdots e_k.$$

With the use of $e_{h'} = e_{h'}\cdots e_k \cdots e_{h'}
= \wh{w_{kh'}} e_k \wh{w_{h'k}}$ this becomes
\begin{eqnarray*}
\wh{i} \wh{w_{kj}} e_k
& = &\wh{w_{h'j}} e_{h'}\wh{i} e_h \cdots e_k 
 = \wh{w_{h'j}}\wh{i}  e_{h'} \wh{w_{kh'}} 
 = \wh{w_{h'j}}  e_{h'}\wh{i} \wh{w_{kh'}} \\
& = &  \wh{w_{h'j}} \wh{w_{kh'}}e_{k}\wh{w_{h'k}} \wh{i}\wh{w_{kh'}} 
 = \wh{w_{kj}} e_{k} \wh{w_{h'k}} \wh{i}\wh{w_{kh'}}.
\end{eqnarray*}
It is easy to verify that $\wh{w_{h'k}} \wh{i}\wh{w_{kh'}}$ commutes
with $e_{k}$.

We are left with the case where $i$ is not in $\gamma$ but is adjacent
to $j$, an end node of
$\gamma$.
Then
$\wh{i} \wh{w_{kj}} e_k
 = \wh{i} \wh{w_{\a_j,k}} e_k 
 = \wh{w_{\a_i+\a_{j},k}} e_k$.
This ends the proof of the equalities involving
$\wh{i} \wh{w_{kj}} e_k$.

\np
We now consider $e_{i} \wh{w_{kj}} e_k$.
If $i=j$, we have trivially
$e_{i} \wh{w_{kj}} e_k \in \wh{w_{kj}} Z_k^{(0)}$.
So let $i \ne j$.
If $i\not\sim j$ we find
$e_{i} \wh{w_{kj}} e_k = e_{i} e_j\cdots \wh{w_{kj}} e_k \in I_2$.
So assume $i\sim j$.

If $i$ occurs in $\gamma$, the path $\gamma$ begins with $j\sim i$ and
so
\begin{eqnarray*}
e_i \wh{w_{kj}} e_k &=& e_ie_je_i\cdots e_k = e_i\cdots e_k
=\wh{w_{ki}}e_k
\end{eqnarray*}
and if $i$ does not occur in $\gamma$, we have
$e_i \wh{w_{kj}} e_k = e_ie_j\cdots e_k = \wh{w_{ki}}e_k$.
\end{proof}

\np
Let $i$ be a node of $M$ and $\beta\in\Phi^+$.
We shall use the following notation.
\begin{itemize}
\item
$\Geod(i,\beta)$ is the set of nodes of the shortest path from
$i$ to a node in the support of $\beta$ that are not in the support
themselves.
So $\Geod(i,\beta) = \emptyset$ if $i\in\Supp(\b)$.
\item $\Proj(i,\beta)$ is the node in the support of $\beta$ nearest $i$.
So $\Proj(i,\beta) = i$ if $i\in\Supp(\b)$.
\item
$C_{\beta,i}$ is the coefficient of $\a_i$ in the expression of
$\b$ as a linear combination of the fundamental roots. So $\beta = \sum_i C_{\beta,i} \a_i$.
\item $J_{\b,k}$ is the subset of $M$ of all nodes $j$ such that $(\a_j,\b)
= 1$ and $\wh{j} \wh{w_{\b-\a_j,h}} = \wh{w_{\b,h}}$, where $h = \Proj(\b,k)$.
This set is empty only if $\b$ is a fundamental root.
\end{itemize}

For $i$ a node of $M$, denote by $i^\perp$ the set of all nodes distinct and
non-adjacent to $i$.

\begin{Lm} \label{finalcase}
Let $\b$ be a root and let $k$ be a node of $M$ such that $i = \Proj(\b,k)$
satisfies $(\a_i,\b) =0$ and $C_{\b,i} = 1$. If $J_{\b,k}\cap i^\perp =
\emptyset$
then
\begin{eqnarray*}
\wh{i}\wh{w_{\b,k}}e_k &=&
\wh{w_{\b,k}}^{-op}e_k\wh{w_{\b,k}}^{op}\wh{i}\wh{w_{\b,k}} \in \wh{w_{\b,k}}^{-op}Z_k^{(0)}.
\end{eqnarray*} 

\end{Lm}

\begin{proof}
We only have to prove that $e_k\wh{w_{\b,k}}^{op}\wh{i}\wh{w_{\b,k}}$ belongs to $Z_k^{(0)}$.
Moreover,
$$e_k\wh{w_{\b,k}}^{op}\wh{i}\wh{w_{\b,k}} = 
e_k\wh{w_{ik}}\wh{w_{\b,i}}^{op}\wh{i}\wh{w_{\b,i}}\wh{w_{ki}}$$
and $J_{\b,k} = J_{\b,i}$, so,
by Proposition \ref{Z_i-prop}(ii), 
it suffices to consider the case where $k=i$.

We prove this by induction on the height of $\b$. The smallest possible root
that satisfies the conditions of the lemma is a root of the form
$\a_j+\a_i+\a_h$ with $j \adj i \adj h$. In this case $\wh{w_{\b,i}} =
\wh{hj}$.
Straightforward computations give
$$
e_i\wh{w_{\b,i}}^{op}\wh{i}\wh{w_{\b,i}}
 = e_i\wh{jh}\wh{i}\wh{hj}
= e_i\wh{jihij}
 =e_i \wh{w_{ji}}\wh{h}\wh{w_{ij}}
 = e_i \wh{w_{ij}}^{op}\wh{h}\wh{w_{ij}},
$$
which belongs to $Z_k^{(0)}$ by definition.

Let $\b$ be a
positive root of height at least 4 and assume that the lemma holds for all
positive roots of height less than $\het(\b)$.  Now $\wh{w_{\b,k}}e_k =
\wh{w_{\b,i}}e_i \cdots e_k$ with no $i$ in $w_{\b,i}$. Let $j \in J_{\b,k}$.
Then, by the hypothesis $J_{\b,k}\cap i^\perp = \emptyset$, we have $i\sim j$.
Clearly $w_{\b,i} = j w_{\b-\a_j,i}$.  As $(\a_i,\b) = 0$ and $C_{\b,i} = 1$,
the sum of $C_{\b,j}$ for $j$ running over the neighbors of $i$ in $M$, must be
$2$. Hence there are either two nodes $j$, $h$ say, in $M$ with $C_{\b,j} =
C_{\b,h} =1$ or there is a single node $j$ of $M$ adjacent to $i$ 
with $C_{\b,j} = 2$. In the former case, as $\het(\b)\ge 4$, there is an end
node $p$ of $\b$ distinct from $j$, $i$, $h$ and non-adjacent to $i$ with
$C_{\b,p} = 1$, which implies $(\a_p,\b) = 1$, whence $p\in J_{\b,i}\cap
i^\perp$, a contradiction. Hence $i$ is an end node of $\b$ and has a
neighbor $j$ with $C_{\b,j}= 2$ and $(\a_j,\b) = 1$. This implies that
$\wh{w_{\b-\a_j,i}} = \wh{w_{\c,j}}\wh{j}$, where $\c=\b-\a_i-\a_j$.  As $(\a_j,\c) =
  0$ and $ J_{\c,j}\cap j^\perp \subseteq J_{\b,i}\cap i^\perp = \emptyset$, we
  can apply induction to find
$e_j\wh{w_{\c,j}}^{op}\wh{j}\wh{w_{\c,j}}$ belongs to $Z_j^{(0)}$.
Consequently,
\begin{eqnarray*}
e_i\wh{w_{\b,i}}^{op}\wh{i}\wh{w_{\b,i}} 
&=&e_i \wh{j} \wh{w_{\c,j}}^{op} \wh{jij} \wh{w_{\c,j}} \wh{j}   
= e_i \wh{j} \wh{w_{\c,j}}^{op} \wh{iji} \wh{w_{\c,j}} \wh{j}   
= e_i\wh{ji} \wh{w_{\c,j}}^{op} \wh{j} \wh{w_{\c,j}} \wh{ij}  \\ 
&\in& \wh{w_{ji}} Z_j^{(0)} \wh{w_{ij}} 
=  Z_i^{(0)}.
\end{eqnarray*} 
\end{proof}

\begin{Lm}\label{aib=0}
Let $\b$ be a root and let $i$ be a node with $(\a_i,\b)=0$. Then the following hold.
\begin{enumerate}[(i)]
\item If $j$ is a node in $J_{\b,k}\cap i^\perp$ then
$\wh{i} \wh{w_{\b,k}}e_k = \wh{j} \wh{i} \wh{w_{\b-\a_j,k}}e_k$.
\item If
$i = \Proj(\b,k)$ and $C_{\b,i} = 1$ and $J_{\b,k}\cap i^\perp=\emptyset$,
then
$\wh{w^{op}_{\b,k}}\wh{i}\wh{w_{\b,k}}\in Z_k^{(0)}$ and
$$\wh{i}\wh{w_{\b,k}}e_k =
\wh{w_{\b,k}}^{-op}e_k(\wh{w^{op}_{\b,k}}\wh{i}\wh{w_{\b,k}}).$$
\item If $i \ne \Proj(\b,k)$ or $C_{\b,i}>1$, then,
for $j\in J_{\b,k}\setminus i^\perp$,
$$\wh{i} \wh{w_{\b,k}}e_k= \wh{ji} \wh{j} \wh{w_{\b-\a_j-\a_i,k}}e_k.$$
\end{enumerate}
In each case, $\wh{i} \wh{w_{\b,k}}e_k\in \wh{w_{\b,k}}Z_k^{(0)}$.
\end{Lm}

\begin{proof}
(i).
Straightforward from $\wh{ij} = \wh{ji}$.

For the remainder of the proof,
we can and will assume there is a node $j$ with
$(\a_j,\b) = 1$, $w_{\b,k} = r_jw_{\b-\a_j,k}$ and $i \adj j$. Then $(\a_i,\b-\a_j)=1$. 

\nl (ii). This follows from Lemma~\ref{finalcase}.

\nl (iii).
Here $\wh{w_{\b,k}} = \wh{ji}\wh{w_{\b-\a_j-\a_i,k}}$ and the statement follows from the
braid relation $\wh{iji} =\wh{jij}$.
\end{proof}

\begin{Thm}\label{doubleI}
Let $B$ be a BMW-algebra of type $M\in \ADE$, let $\b\in \Phi^+$, and
let $i$, $k$ be nodes of $M$.

If $(\a_i,\b) = -1$, then
$$
\wh{i}\wh{w_{\b,k}}e_k = 
\left\{
\begin{tabular}{rcl}
$\wh{w_{\b+\a_i,k}}e_k$&if&$i\not\in \Geod(k,\beta)$,\\
$\wh{w_{\b+\a_i,k}}e_k-m\wh{w_{\b,k}}e_k+m\wh{w_{\a_i,k}}e_k\wh{w_{\beta,h}}
$&if&$i\in \Geod(k,\beta)$ and\\
&&$h = \Proj(k,\beta)$.
\end{tabular}
\right.
$$

If $(\a_i,\b) = 1$, then
$$
\wh{i}\wh{w_{\b,k}}e_k = 
\left\{
\begin{tabular}{rcl}
$\wh{w_{\b-\a_i,k}}e_k-m\wh{w_{\b,k}}e_k+ml^{-1}e_i\wh{w_{\b-\a_i,k}}e_k$&if
&$i\in J_{\b,k}$,\\
$\wh{w_{\b-\a_i,k}}e_k$&if&$i\not\in J_{\b,k}$.
\end{tabular}
\right.
$$

If $(\a_i,\b) = 0$, then
$$
\wh{i}\wh{w_{\b,k}}e_k = 
\left\{
\begin{tabular}{rcl}
$\wh{w_{\b,k}}e_k (\wh{w_{\b,k}}^{-1} \wh{i} \wh{w_{\b,k}})$&if&$i\not\in \Supp(\b)$,\\
$\wh{ji}\wh{w_{\b-\a_j,k}}e_k$&if&$j\in J_{\b,k}\cap i^\perp$,\\
$\wh{w_{\b,k}}^{-op}e_k(\wh{w^{op}_{\b,k}}\wh{i}\wh{w_{\b,k}})$&if&$C_{\b,i} =
1$, $i = \Proj(\b,k)$,\\
&& and $J_{\b,k}\cap i^\perp=\emptyset$,\\ 
$\wh{jij} \wh{w_{\b-\a_j-\a_i,k}}e_k$&if&$j\in J_{\b,k}\setminus i^\perp$ and $i \in J_{\b-\a_j,k}$.
\end{tabular}
\right.
$$

If $(\a_i,\b) = 2$, then $\b = \a_i$ and
$
\wh{i}\wh{w_{\b,k}}e_k = l^{-1}\wh{w_{\b,k}}e_k$.

In each case, the result is in $\wh{w_{\c,k}}Z_k^{(0)} +
m\wh{w_{\b}}Z_k^{(0)}
+ m\wh{w_{\a_i,k}}Z_k^{(0)} + I_2$,
where $\c = \b$ if $\b=\a_i$ and 
$\c = r_i\b$ otherwise.
\end{Thm}

\begin{proof}
By Lemma~\ref{H0} the theorem holds for all fundamental roots $\b$ in $\Phi^+$.
Suppose $\b$ is a non-fundamental root in $\Phi^+$, and
consider $\wh{i}\wh{w_{\b,k}}e_k$.
Now $(\a_i,\b)<2$, for otherwise
$\b=\a_i$.
First let $(\a_i,\b)=1$. If $i\in J_{\b,k}$, then 
$$
\wh{i}\wh{w_{\b,k}} e_k = \wh{i}^2\wh{w_{\b-\a_i,k}} e_k =
\wh{w_{\b-\a_i,k}}e_k - m \wh{w_{\b,k}} e_k +ml^{-1}
e_i\wh{w_{\b-\a_i,k}} e_k.$$ 

Assume $i\not\in J_{\b,k}$ then $i = \Proj(k,\b)$ and $C_{\b,i}=1$.
There must be a single node $j\in \Supp(\b)\setminus i^\perp$ with $C_{\b,j} = 1$,
and the remaining nodes in the support of $\b$ are on the side of $j$
in $M$ other 
than $i$. This means $\wh{w_{\b,k}} = \wh{u}\wh{j} \wh{w_{\a_i,k}} $ where
the elements in $u$ are on the side of $j$ other than $i$ and so $i$
commutes with $u$. Now $\wh{i}\wh{u}\wh{j}\wh{w_{\a_i,k}} =
\wh{u}\wh{i}\wh{j}\wh{w_{\a_i,k}} =\wh{u}\wh{w_{\a_j,k}} $ so
$\wh{i}\wh{w_{\b,k}}e_k = \wh{w_{\b-\a_i,k}}e_k $ as required. 

\np
Next let $(\a_i,\b)=0$ and assume $i$ is not in the support of $\b$.
Put $h = \Proj(k,\b)$ and $\rho = \Geod(k,\b)$.
If $i$ is not in $\rho$ and not adjacent to an element of $\rho$, then
$\wh{i}$ commutes with $\wh{w_{\b,k}}$ so $\wh{w_{\b,k}}^{-1}
\wh{i} \wh{w_{\b,k}}= \wh{i}$ and
$\wh{i} \wh{w_{\b,k}}e_k = \wh{w_{\b,k}} e_k\wh{i}$.

If $i$ is in $\rho$ or adjacent to an element of $\rho$, then $\wh{i}$
commutes with $\wh{w_{\b,h}}$ where $\wh{w_{\b,k}} =
\wh{w_{\b,h}}\wh{w_{\a_h,k}}$. Now
$\wh{w_{\b,k}}^{-1} \wh{i} \wh{w_{\b,k}} = \wh{w_{\a_h,k}}^{-1} \wh{i}
\wh{w_{\a_h,k}}$.
We know that $i \not\adj h$ so $\wh{w_{\a_h,k}}^{-1}\wh{i} \wh{w_{\a_h,k}}e_k
\in Z_k^{(0)}$ by Lemma~\ref{H0}. We conclude $\wh{i} \wh{w_{\b,k}}e_k
= \wh{w_{\b,k}} \wh{w_{\b,k}}^{-1}\wh{i} \wh{w_{\b,k}}e_k
= \wh{w_{\b,k}} \wh{w_{\a_h,k}}^{-1}\wh{i} \wh{w_{\a_h,k}}e_k 
= \wh{w_{\b,k}} e_k( \wh{w_{\a_h,k}}^{-1}\wh{i} \wh{w_{\a_h,k}} )
\in \wh{w_{\b,k}} Z_k^{(0)}$.

\np If $(\a_i,\b)=0$ with $i\in \Supp(\b)$, then the assertion follows from Lemma \ref{aib=0}.

\np 
Finally let $(\a_i,\b)=-1$. Here $\wh{i} \wh{w_{\b,k}}e_k =
\wh{w_{\b+\a_i,k}}e_k$ by definition if $i$ is not in $\Geod(k,\b)$.
So suppose $i\in \Geod(k,\beta)$. 
Write
$h = \Proj(k,\beta)$. Since
$(\a_i,\beta) = -1$, we must have $i\sim h$.
Therefore $\wh{w_{\b,k}} = \wh{w_{\beta,h}}\wh{w_{kh}}$ and
$\wh{w_{\b,k}}e_k = \wh{w_{\beta,h}}e_he_i\cdots e_k$.
The set $\Supp(\b) \setminus \{h\}$ is a connected component of the Dynkin diagram
connected to $h$ and disconnected from 
$ \Geod(k,\beta)$.
Hence $\wh{h}$ does not appear in $\wh{w_{\beta,h}}$.
This means $\wh{i}$ commutes
with $\wh{w_{\beta,h}}$.
Moreover, by definition of $w_{\b,h}$, we have
$ \wh{w_{\beta,h}}\wh{h} = 
\wh{w_{\b+\a_i,i}}$
and so
$ \wh{w_{\beta,h}}\wh{h}\wh{w_{ki}} = 
\wh{w_{\b+\a_i,k}}$.
Consequently,
by (\ref{7-eq}),
\begin{eqnarray*}
\wh{i}\wh{w_{\beta,h}}e_he_i\cdots e_k &=&
\wh{w_{\beta,h}}\wh{i}e_he_i\cdots e_k 
= \wh{w_{\beta,h}}(\wh{h} + m(1-e_h))e_i\cdots e_k \\ 
&=& \wh{w_{\beta,h}}\wh{h}\wh{w_{ki}}e_k
+
m\wh{w_{ki}} e_k \wh{w_{\beta,h}}
-m\wh{w_{\beta,h}}e_h\wh{kh} \\ 
&=&\wh{w_{\b+\a_i,k}}e_k
+m\wh{w_{k\,i}}e_k\wh{w_{\beta,h}}
-m \wh{w_{\b,k}}e_k
\end{eqnarray*} 

\end{proof}

\begin{Cor}\label{doublee}
Let $B$ be a BMW-algebra of type $M\in \ADE$, let $\b\in \Phi^+$, and
let $i$, $k$ be nodes of $M$.
\begin{enumerate}[(i)]
\item \label{it1}
$\wh{w_{\b,k}}^{-op}e_k  \in  \wh{w_{\b,k}}e_k + m \sum_{\het(\c) < \het(\b)}
\wh{w_{\c,k}}Z_k^{(0)} +I_2$,
\item \label{it2}
$e_i \wh{w_{\b,k}} e_k  \in \wh{w_{\a_i,k}} Z_k^{(0)} + I_2$,
\item \label{it3}
$\wh{i} \wh{w_{\b,k}} e_k  \in \sum_{\c \in H_{\b,i}} \wh{w_{\c,k}}  Z_k^{(0)}
+ I_2$,
\end{enumerate}
where $H_{\b,i} = \{ \delta \in \Phi^+ \mid \het(\delta) < \het(\b) \} \cup \{
\b,\b+\a_i \}\cap\Phi^+$.
\end{Cor}

\begin{proof}
We prove the statements simultaneously by induction on the height of $\b$.
If $\b$ is a fundamental root then statement (i) holds by
Lemma~\ref{ij-paths} and the statements (ii) and (iii) by Lemma~\ref{H0}. 

Let $\b \in \Phi^+$ with $\het(\b) \ge 2$ and assume the lemma holds for
all $\c \in \Phi^+$ with $\het(\c) < \het(\b)$.
Let $i,k$ be nodes and consider $\wh{w_{\b,k}}^{-op}e_k$, $e_i \wh{w_{\b,k}}
e_k$ and $\wh{i} \wh{w_{\b,k}} e_k$.
There is (at least one) $j$ such that $\wh{w_{\b,k}} =
\wh{j}\wh{w_{\b-\a_j,k}}$; then $\het(\b-\a_j) 
= \het(\b) - 1$. 
Now 
\begin{eqnarray*}
\wh{w_{\b,k}}^{-op}e_k &=& \wh{j}^{-1}\wh{w_{\b-\a_j,k}}^{-op}e_k \\
&\in& (\wh{j}+ m -me_j) (\wh{w_{\b-\a_j,k}}e_k + m \sum_{\het(\c) <
\het(\b-\a_j)}
\wh{w_{\c,k}}Z_k^{(0)} +I_2) \\
&=& \wh{w_{\b,k}}e_k + m\wh{w_{\b-\a_j,k}}e_k - me_j\wh{w_{\b-\a_j,k}}e_k +
m \sum_{\het(\c) < \het(\b-\a_j)} \wh{j}\wh{w_{\c,k}}Z_k^{(0)} \\
&& + m^2 \sum_{\het(\c) < \het(\b-\a_j)} \wh{w_{\c,k}}Z_k^{(0)}  - m^2 \sum_{\het(\c)
< \het(\b-\a_j)}
e_j\wh{w_{\c,k}}Z_k^{(0)} + I_2 \\
&\subseteq & \wh{w_{\b,k}}e_k + m \sum_{\het(\c) < \het(\b)}\wh{w_{\c,k}}Z_k^{(0)}
+ m \sum_{\het(\c) < \het(\b)}
e_j\wh{w_{\c,k}}Z_k^{(0)} + I_2 \\
&\subseteq & \wh{w_{\b,k}}e_k + m \sum_{\het(\c) < \het(\b)}\wh{w_{\c,k}}Z_k^{(0)} + I_2.
\end{eqnarray*}
To see that $\sum_{\het(\c) < \het(\b-\a_j)}\wh{j}\wh{w_{\c,k}}Z_k^{(0)}$ is
contained in $\sum_{\het(\c) < \het(\b)} \wh{w_{\c,k}}Z_k^{(0)}$, observe that
by the induction hypothesis on (iii) we have $\wh{j}\wh{w_{\c,k}}e_k \in \sum_{\delta \in
  H_{\c,i}} \wh{w_{\delta,k}} Z_k^{(0)} + I_2$. Here $\het(\delta) \le
\het(\c)+1 < \het(\b)$ while $\het(\c) < \het(\b-\a_j)$.  The sum
$\sum_{\het(\c) < \het(\b)} e_j\wh{w_{\c,k}}Z_k^{(0)}$ is in
$\wh{w_{\a_j,k}}Z_k^{(0)}$ by our induction hypothesis on (ii)
and this gives (i) for $\b$.

Now focus on $e_i \wh{w_{\b,k}} e_k = e_i \wh{j}\wh{w_{\b-\a_j,k}} e_k$.
If $i =j$ then, by the induction hypothesis, $e_i \wh{w_{\b,k}} e_k = l^{-1}
e_i 
\wh{w_{\b-\a_j,k}} e_k \in \wh{w_{\a_i,k}} Z_k^{(0)} + I_2$.
If $i \not\adj j$ then $e_i \wh{j} \wh{w_{\b-\a_j,k}} e_k = \wh{j} 
e_i \wh{w_{\b-\a_j,k}} e_k \in \wh{j}\wh{w_{\a_i,k}} Z_k^{(0)} + I_2$ 
and by Lemma~\ref{H0} this is contained in $\wh{w_{\a_i,k}} Z_k^{(0)} + I_2$.

So, for the remainder of the proof, we may (and shall) assume $i \adj j$.  By
(\ref{7-eq}), we have $e_i \wh{j} = e_i e_j \wh{i} + m e_i e_j - me_i$, so
$$e_i\wh{w_{\b,k}}e_k=e_i\wh{j}\wh{w_{\b-\a_j,k}}e_k= e_i e_j \wh{i}\wh{w_{\b-\a_j,k}} e_k + m e_i e_j\wh{w_{\b-\a_j,k}} e_k -
me_i\wh{w_{\b-\a_j,k}} e_k.$$
By our induction hypothesis the last two terms are in $\wh{w_{\a_i,k}} Z_k^{(0)} + I_2$.
This leaves the first term, $e_i e_j \wh{i}\wh{w_{\b-\a_j,k}} e_k$.
Because $(\b-\a_j,\a_i) = (\b,\a_i)+1$ the inner product of $\a_i$ with
$\b-\a_j$ can only take values $0$, $1$, and $2$ and thus $H_{\b-\a_j,i}$
consists of roots with height at most $\het(\b-\a_j)$.

The induction hypothesis on (iii) now gives
$\wh{i}\wh{w_{\b-\a_j,k}} e_k \in \sum_{\c \in
H_{\b-\a_j,i}}\wh{w_{\c,k}}Z_k^{(0)}$ where $\het(\c) < \het(\b)$ for all $\c$.
By applying the induction hypothesis twice we obtain
$$e_i e_j \wh{i}\wh{w_{\b-\a_j,k}} e_k \in \sum_{\het(\c) \in
H_{\b-\a_j,i}}e_i e_j\wh{w_{\c,k}}Z_k^{(0)} + I_2 \subseteq e_i \wh{w_{\a_j,k}}Z_k^{(0)}
\subseteq \wh{w_{\a_i,k}}Z_k^{(0)} + I_2.$$
This establishes (\ref{it2}).
Finally consider $\wh{i} \wh{w_{\b,k}}e_k$.
If $(\a_i,\b) = -1$ then $\b+\a_i \in \Phi^+$ and the statement holds by
Theorem~\ref{doubleI}.
Also, if $(\a_i,\b) = 1$ then Theorem~\ref{doubleI}
applies. Here $e_i\wh{w_{\b-\a_i,k}}e_k \in \wh{w_{\a_i,k}}Z_k^{(0)} + I_2$ by
the induction hypothesis for (\ref{it2}).

For the remainder of the proof we assume $(\a_i,\b) = 0$. Again
Theorem~\ref{doubleI} gives
an expression for $\wh{i}\wh{w_{\b,k}}e_k $ in each of the four cases discerned.
In the first cases, where $i\not\in \Supp(\b)$, the statement is immediate from this expression. By our induction hypothesis
for (iii) the second case gives an expression contained in $\sum_{\c \in
H_{\b-\a_j,i}} \wh{j}\wh{w_{\c,k}} Z_k + I_2$ whence in $\sum_{\c \in H_{\b,i}}
\wh{w_{\c,k}} Z_k + I_2$. Now the fourth case goes by the same argument
and only the third case remains to be verified. Above we have shown that
$\wh{w_{\b,k}}^{-op}e_k \in
\wh{w_{\b,k}}e_k + m \sum_{\het(\c) < \het(\b)}\wh{w_{\c,k}}Z_k^{(0)} +I_2$ and
that completes the proof.
\end{proof}

We shall use the following lemma to derive an upper bound for $\dim (Z_i)$ from
Theorem
\ref{doubleI}.

\begin{Lm}\label{PID-lm}
Suppose $F$ is a field, $E$ is a subring of $F$ which is a principal ideal
domain. If $V$ is a vector
space over $F$ and $V^{(0)}$ is an $E$-submodule of $V$ containing a spanning set of
$V$, then $V^{(0)}$ is a free $E$-module on a basis of $V$.
Moreover, if $a\in E$ generates a maximal ideal of $E$, then
$$\dim_{F}(V) = \dim_{E/aE}(V^{(0)}/aV^{(0)}).$$
\end{Lm}

\begin{proof}
  As $E$ is a principal ideal domain, it is well known, see \cite{DF}, Theorem
  12.5, that each $E$-module of finite rank without torsion is free. Applying
  this observation to $V^{(0)}$, we let $X$ be a basis of the $E$-module
  $V^{(0)}$.  By the hypothesis that $V^{(0)}$ spans $V$, it is also a basis
  of $V$, so $\dim_F (V) = |X|$.  On the other hand, $X$ maps onto a basis of
  $V^{(0)}/aV^{(0)}$ over $E/aE$ (for, it clearly maps onto a spanning set and
  if $\sum_{x\in X}\lambda_x x = 0\mod aV^{(0)}$ for $\lambda_x\in E$, then,
  as $V^{(0)} = E\, X$, with $X$ a basis, we have $\lambda_x = 0\mod a$ for each $x\in X$, so
  the linear relation in $V^{(0)}/aV^{(0)}$ is the trivial one).  This proves $\dim_F(V)
  = |X| = \dim_{E/aE}(V^{(0)}/aV^{(0)})$.
\end{proof}

\begin{Cor}\label{BmodI2spanning}
Let $M\in \ADE$
and let $i$ be a node of $M$. Then
$\wh{D_i}Z_i\wh{D_i}^{op}$ is a linear spanning set for $I_1/I_2$.
Moreover, the dimension of $Z_i$ is at most $|W_C|$.
\end{Cor}

\begin{proof}
By Lemma \ref{conj-induction} $I_1$ is spanned by a set of multiples of $e_i$
by generators $g_j$, so $I_1 = Be_iB$.
According to Theorem \ref{doubleI} and Corollary \ref{doublee}, 
$Be_i = \wh{D_i}Z_i+I_2$. Applying Remark \ref{opLm}, we derive
from this that $e_iB = Z_i\wh{D_i}^{op}+I_2$
(observe that $Z_i$ and $I_2$ are invariant under the anti-involution).
Therefore, $I_1 = Be_iB = \wh{D_i}Z_i\wh{D_i}^{op}+I_2$.

It remains to establish that the dimension of $Z_i\mod I_2$ is at most $|W_C|$.
To this end we consider the integral versions $Z_i^{(0)}$  and $B^{(0)}$ of
$Z_i$ and $B$ over $E = \Q(x)[l^{\pm}]$ defined at the beginning of
Section \ref{sec:I1modI2},
and look at the quotients modulo $(l-1)$.
Observe that, by (\ref{def-x-eq}), $m$ belongs to the ideal $(l-1)E$.

A careful inspection of the identities in Theorem \ref{doubleI} and Corollary
\ref{doublee}, shows $B^{(0)}e_i = \wh{D_i}Z_i^{(0)} + I_2$, and $e_i
B^{(0)}e_i = Z_i^{(0)}+I_2$.  Since $Be_i$ is linearly spanned by the set $
\wh{D_i}Z_i^{(0)}
$ mod $I_2$, it is linearly spanned by
$B_i^{(0)}e_i$ mod $I_2$.  Consequently, $Z_i = e_iBe_i$ is linearly spanned
by $e_iB_i^{(0)}e_i +I_2$, whence by $Z_i^{(0)}$ mod $I_2$.

For brevity of notation, we set $m_1 = l-1$. (The remainder of the proof would
also work for $m_1 = l+1$.)
Since $x^{-1}e_i$ is a central idempotent belonging to $Z_i^{(0)}$, we have
$m_1 B^{(0)}\cap (Z_i^{(0)} +I_2) = m_1 e_iB^{(0)}e_i \cap (Z_i^{(0)}
+I_2) = m_1(Z_i^{(0)}+I_2) \cap (Z_i^{(0)} + I_2) = m_1 (Z_i^{(0)} + I_2)
$.  Therefore, the quotient $Z_i^{(0)}/m_1Z_i^{(0)}$, viewed as a vector
space over $\Q(x)$, is isomorphic to $ (Z_i^{(0)} + m_1
B_i^{(0)} +I_2)/(m_1B_i^{(0)} + I_2)$.  But this algebra 
is readily seen to be a quotient of a
subalgebra of the group algebra over $\Q(x)$ of the stabilizer in $W$ of the
simple root $\a_i$,  for the image of $\{ \wh{w}\mid w\in W\}$ modulo $m_1B^{(0)}$
is the group $W$ and
the image of the algebra $Z_i^{(0)}$ is generated by the products of the elements of
the form $w_{ji}r_{k}w_{ji}$ for $j$ and $k$ distinct non-adjacent nodes of
$M$,
all of which are contained in the stabilizer in $W$ of $\a_i$.
Consequently, the dimension of $Z_i^{(0)}/m_1 Z_i^{(0)}$ over $\Q(x)$ is at
most $|W_C|$, the order of the stabilizer in $W$ of $\a_0$ (a group conjugate
to the stabilizer in $W$ of $\a_i$).
By Lemma \ref{PID-lm}, applied with $F = \Q(x,l)$, $E = \Q(x)[l^{\pm}]$,
$V = Z_i$, $V^{(0)} = Z_i^{(0)}$, and $a = m_1$,
we see that
$Z_i $ has dimension at most $|W_C|$ over $\Q(l,x)$.
\end{proof}

\section{Generalized Lawrence-Krammer representations} \label{sec:krammer}
In this section we construct the analog of the Lawrence-Krammer representation
of $A$ with coefficients in $Z_0$, the Hecke
algebra of type $C$, where $C$ is the parabolic of the highest root
centralizer.  We show the representation factors through $B/I_2$.  By taking
an irreducible representation of $Z_0$, we find an irreducible representation of
$B/I_2$. Finally, by counting dimensions of irreducible representations, we
are able to conclude that all representations of $B/I_2$ that do not vanish on
$I_1$ are of this generalized Lawrence-Krammer type, and we can finish the
proof of Theorem \ref{main-thm}.

Since the construction for disconnected $M$ is a direct sum of the
representations of $B$ for the distinct connected components, we simply take
$M$ to be connected, so $M\in\ADE$.
We let $\Phi$ be the root system in $\R^n$ of type $M$, and
denote by $\a_1,\ldots,\a_n$ the fundamental roots
corresponding to the reflections $r_1,\ldots,r_n$,
respectively. As usual, by $\Phi^+$ we denote the set of positive roots
in $\Phi$.

For a root $\b$, the set of roots $\{\gamma\in\Phi\mid (\b,\gamma) = 0\}$ is
also a root system. Its type can be read off from $M$ as follows: the extended
Dynkin diagram $\widetilde K$ of the connected component $K$ of $M$ involving
$\b$ (i.e., having nodes in the support of $\beta $) has a single node
$\alpha_0$ in addition to those of $K$; now take $C$ to consist of all nodes
of $M$ that are not connected to $\a_0$.  Then the type of the roots
orthogonal to $\b$ is $M|_C$. In fact, if $\b=\a_0$, then $\{\a_i\mid i\in
C\}$ is a set of fundamental roots of the root system $\{\gamma\in\Phi\mid
(\b,\gamma) = 0\}$.
For $\A_n$ with $\b=\a_0$ this is the diagram of type $\A_{n-2}$ on $\{2,\ldots,n-1\}$,
for $\D_n$, it is the diagram of type $\A_1\D_{n-2}$ on
$\{1\}\cup\{3,\ldots,n\}$,
for $\E_6$ it is the diagram of type $\A_5$ on $\{1,3,4,5,6\}$,
for $\E_7$ it is the diagram of type $\D_6$ on $\{2,3,4,5,6,7\}$,
and for $\E_8$ it is the diagram of type $\E_7$ on $\{1,2,3,4,5,6,7\}$.

Recall the coefficients of $Z_0$ are in $\Q(l,x)$.  We take the coefficients of our representation in
the Hecke algebra $Z_0^{(0)}$ of type $M|_C$ over the subdomain
$\Q[l^{\pm1},m]$ of $\Q(l,x)$, where $m$ is defined in (\ref{def-x-eq}). Observe
that the fraction field of $\Q(l,m)$ coincides with $\Q(l,x)$. The
generators $z_i$ $(i\in C)$ of $Z_0^{(0)}$ satisfy the
quadratic relations $z_i^2+mz_i-1 = 0$.  For the proof of irreducibility at
the end of this section, we need however a smaller version of this Hecke
algebra, namely the subalgebra $Z_0^{(1)}$ with same generators $z_i$, but
over $\Q[m]$. Thus, $Z_0^{(0)} =Z_0^{(1)} \Q[l^{\pm1}] $.

By Lemma \ref{Hbeta-lm}, the element $h_{\b,i} $ of $A$ defined in (\ref{hbetai-df}),
where $\b\in\Phi^+$ and $i$ is a node with $(\a_i,\b) =
0$, maps onto an
element of $Z_0^{(1)}$ upon substitution of $s_j$ by $z_j$ and $s_j^{-1}$ by
$z_j+m$.  We shall also write $h_{\b,i}$ for the image of this element in
$Z_0^{(1)}$.

We write $V^{(0)}$ for the free right $Z_0^{(0)}$ module with basis $x_\b$
indexed by $\b\in\Phi^+$. The connection with \cite{CW} is given by $m = r -
r^{-1},\quad l= 1/(tr^3)$.
Recall that 
$A^+$ is the positive monoid of $A$.

\begin{Thm} \label{adethmgen2}
Let $M\in \ADE$ and let $A$ be the Artin group of type $M$.
Then, for each $i\in\{1,\ldots,n\}$
and each $\b\in\Phi^+$, there are elements $T_{i,\b}$
in $Z_0^{(1)}$
such that
the following map on the generators of $A$ determines a representation
of $A$ on $V^{(0)}$. 
$$
s_i\mapsto \sigma_i = \tau_i+l^{-1}T_i,
$$
where $\tau_i$ is determined by
\begin{eqnarray*}
\tau_i(x_\b )&=&
\left\{
\begin{tabular}{lr}
$0$ & \mbox{ if }\quad $\phantom{-}(\a_i,\b)=2$\phantom{,}\\
$ x_{\b-\a_i} $ & \mbox{ if }\quad $\phantom{-}(\a_i,\b)=1$\phantom{,}\\
$ x_{\b}h_{\b,i}$ & \mbox{ if }\quad $\phantom{-}(\a_i,\b)=0$\phantom{,}\\
$ x_{\b+\a_i} -m x_{\b} $ &\phantom{-} \mbox{ if }\quad $(\a_i,\b)=-1$,\\
\end{tabular}
\right.
\end{eqnarray*}
and where
$T_i$ is the
$Z_0^{(0)}$-linear map on $V^{(0)}$ determined by $T_ix_\b = x_{\a_i}T_{i,\b}$
on the generators of $V^{(0)}$ and by $T_{i,\a_i} = 1$. 

When tensored with $\Q(x,l)$, the representation of $A$ on $V^{(0)}$ becomes a
representation on the vector space $V$ which factors through the quotient
$B/I_2$ of the BMW algebra $B$ of type $M$ over $\Q(x,l)$.
\end{Thm}

\np
Throughout this section we use several properties of the elements
$h_{\b,i}$ listed in Lemma \ref{H-lem}. In addition, we shall use
the Hecke algebra relation for the image of $h_{\b,i}$ in $Z_0^{(0)}$:
\begin{eqnarray}
h^{-1}_{\b,i} & = & h_{\b,i} + m. \label{H-rel}
\end{eqnarray}

\np
The proof of the theorem follows the lines of the proof in \cite{CW}.
We shall first describe the part modulo $l^{-1}$ of the representation of the
Artin monoid $A^+$
on $V^{(0)}$.

\begin{Lm}\label{monmorphismlm}
There is a monoid homomorphism $A^+ \to {\rm End}(V^{(0)})$
determined by $s_i\mapsto\tau_i$ $(i=1,\ldots,n)$.
\end{Lm}

\begin{proof}
We must show that, if $i$ and $j$ are not adjacent, then $\tau_i
\tau_j=\tau_j \tau_i$ and, if they are adjacent, then
$\tau_i\tau_j\tau_i=\tau_j\tau_i\tau_j$. We evaluate the 
expressions
on each $x_\b $ and show they are equal. We begin with the case where
$\b=\a_i$. 
Suppose
first that $i$ and $j$ are not adjacent. Then $\tau_i x_{\a_i} = 0$ 
and
$\tau_jx_{\a_i} = x_{\a_i}h_{\b,j}$. Now $\tau_j\tau_ix_{\a_i} 
=0$ and
$\tau_i\tau_j x_{\a_i} = \tau_ix_{\a_i}h_{\b,j}=0$, so the 
result holds.
Suppose next that $i$ and $j$ are adjacent. Then
$\tau_ix_{\a_i} =\tau_jx_{\a_j} =0$ and
$\tau_jx_{\a_i} =-mx_{\a_i} + x_{\a_i+\a_j} $. Now
\begin{eqnarray*}
\tau_i\tau_j\tau_ix_{\a_i}& = &\tau_i\tau_j(0)=0
,\ \ \mbox{ \rm and }
\\
\tau_j\tau_i\tau_j x_{\a_i} &=& \tau_j\tau_i(-mx_{\a_i} + 
x_{\a_i+\a_j} ) 
= \tau_j(x_{\a_i+\a_j-\a_i} ) = \tau_jx_{\a_j} 
= 0.
\end{eqnarray*}
This ends the verification for the case where $\b=\a_i$.
We now divide the verifications into the various cases depending on the inner
products $(\a_i,\b)$ and $(\a_j,\b)$.  By the above, we may assume $(\a_i,\b),
(\a_j,\b)\ne 2$.

\np
First assume that $(\a_i,\a_j)=0$. 
The computations verifying $\tau_i\tau_j=\tau_j\tau_i$ are
summarized in the following table.
The last column indicates the formulas that are used.
\medskip
\begin{center}
\begin{tabular}{|r|r|c|c|}
\hline
$(\a_i,\b)$&$(\a_j,\b)$&$\tau_i\tau_j x_\b=\tau_j\tau_ix_\b$&ref.\\
\hline
$1$\ \ &$1$\ \ & $x_{\b-\a_i-\a_j}$ &\\
$1$\ \ &$-1$\ \ &$x_{\b+\a_j-\a_i}-mx_{\b-\a_i}$&\\
$1$\ \ &$0$\ \ &$x_{\b-\a_i}h_{\b,j}$ &(\ref{H-orth-rel})\\
$0$\ \ &$0$\ \ &$x_\b h_{\b,i}h_{\b,j}$&(\ref{H-com-rel})\\
$0$\ \ &$-1$\ \ &$x_{\b+\a_j}h_{\b,i}-mx_{\b}h_{\b,i}$ &(\ref{H-orth-rel})\\
$-1$\ \ &$-1$\ \ 
&$m^2x_\b-m(x_{\b+\a_i}+x_{\b+\a_j})+x_{\b+\a_i+\a_j}$ &\\
\hline
\end{tabular}
\end{center}

\nl
We demonstrate how to derive these expressions by checking the third
line. 
\begin{eqnarray*}
\tau_i\tau_jx_{\b} &=& \tau_i( x_\b h_{\b,j}) 
= x_{\b-\a_i}h_{\b,j}. 
\end{eqnarray*}
In the other order,
\begin{eqnarray*}
\tau_j\tau_ix_{\b} &=& \tau_j(x_{\b-\a_i}) 
= x_{\b-\a_i}h_{\b-\a_i,j}. 
\end{eqnarray*}
Equality between $h_{\b,j}$ and $h_{\b-\a_i,j}$ follows from
(\ref{H-orth-rel}).

Suppose next that $i\sim j$. The same situation occurs except 
the 
computations are sometimes longer and one case does not occur. This is 
the case where $(\a_i,\b)=(\a_j,\b)=-1$. For then
$\b+\a_i$ is also a root, and $(\b+\a_i,\a_j)=-1-1=-2$. This means 
$\b+\a_i=-\a_j$ and $\b$ is not a positive root.
The table is as follows.

\begin{center}
\begin{tabular}{|r|r|c|c|}
\hline
$(\a_i,\b) $&$(\a_j,\b)$&$\tau_i\tau_j 
\tau_ix_\b=\tau_j\tau_i\tau_jx_\b$&ref.\\
\hline
$1$\ \ &$1$\ \ & $0$ &\\
$1$\ \ &$-1$\ \ &$x_\b h_{\b-\a_i,j}-mx_{\b-\a_i}h_{\b-\a_i,j}$&(\ref{H-pm-rel}) \\
$1$\ \ &$0$\ \ &$x_{\b-\a_i-\a_j} h_{\b,j}$&(\ref{H-10-rel}) \\
$0$\ \ &$0$\ \ &$x_\b h_{\b,i}h_{\b,j}h_{\b,i}$&(\ref{H-braid-rel})\\
$0$\ \ &$-1$\ \ &$x_{\b+\a_i+\a_j}h_{\b,i}-m 
x_{\b+\a_j}h_{\b,i}-mx_{\b}h_{\b,i}^2$ &(\ref{H-0-1-rel}), (\ref{H-rel})\\
$-1$\ \ &$-1$\ \ &{\rm does not occur} &\\
\hline
\end{tabular}
\end{center}
\end{proof}

\np We next study the possibilities for the parameters $T_{k,\b}$ occurring in
Theorem \ref{adethmgen2}.  Recall that there we defined
$\s_k=\tau_k+l^{-1}T_k$, where $T_kx_\b=x_{\a_k}T_{k,\b}$.  We shall introduce
$T_{k,\b}$ as elements of the Hecke algebra
$Z_0^{(0)}$ of type $M|_C$.

\begin{Prop}\label{finalTiRelations}
Set $T_{i,\a_i}=1$ for all $i\in\{1,\ldots,n\}$.
For $\sigma_i\mapsto \tau_i+l^{-1}T_i$ to define a linear 
representation
of the group $A$ on $V$,
it is necessary and sufficient that the equations in Table
\ref{algtable} are satisfied
for each $k,j=1,\ldots,n$ and each $\b\in\Phi^+$.
\end{Prop}

\begin{proof}
The $\s_k$ should satisfy the relations (B1), (B2).
Substituting $\tau_k+l^{-1}T_k$ for $\sigma_k$, we find relations for the
coefficients of $l^{-i}$ with $i=0,1,2,3$. The constant part 
involves
only the $\tau_k$. It follows from Lemma \ref{monmorphismlm} that
these equations are satisfied. We shall derive all of the equations
of Table \ref{algtable} except for (\ref{knotleq}) from the 
$l^{-1}$-linear
part and the remaining one from the $l^{-1}$-quadratic part of the
relations.

The coefficients of $l^{-1}$ lead to
\begin{eqnarray}\label{genTrepeqcom}
T_i\tau_j = \tau_j T_i \ \ \mbox{ and }\ \ \ 
T_j\tau_i = \tau_iT_j \ \ &\mbox{if}&\ \ 
i\not\sim j,\\
\ \ \ \tau_jT_i \tau_j + T_j\tau_i \tau_j
+\tau_j\tau_iT_j
=\tau_iT_j \tau_i
+ T_i\tau_j \tau_i + \tau_i\tau_jT_i
\ \ &\label{genTrepeqadj}\mbox{if}&\ \ 
i\sim j.
\end{eqnarray}

We focus on the consequences of these equations for the $T_{i,\b}$.
First consider the case where $i\not\sim j$. Then 
$\tau_i x_{\a_j} = x_{\a_j} h_{\a_j,i}$ and so, 
for the various values of $(\a_i,\b)$ we find
the following equations
\begin{center}
\begin{tabular}{|r|c|c|}
\hline
$(\a_i,\b)$&
$T_j\tau_i x_{\b} = \tau_iT_j x_{\b}$
&equation\\
\hline
$0$\ \ &
$x_{\a_j} T_{j,\b}h_{\b,i} = 
x_{\a_j} h_{\a_j,i} T_{j,\b}$& $ T_{j,\b}h_{\b,i} = 
h_{\a_j,i} T_{j,\b}$ \\
$1$\ \ &
$x_{\a_j}T_{j,\b-\a_i}=x_{\a_j}h_{\a_j,i}T_{j,\b}$
&$T_{j,\b-\a_i}=h_{\a_j,i}T_{j,\b}$\\
$-1$\ \ &
$x_{\a_j}T_{j,\b+\a_i}-mx_{\a_j}T_{j,\b} = x_{\a_j}h_{\a_j,i}T_{j,\b}$&
$T_{j,\b+\a_i}=h_{\a_j,i}^{-1}T_{j,\b}$\\
$2$\ \ &
$0=x_{\a_j}h_{\a_j,i}T_{j,\b}$ & $0=T_{j,\b}$\\
\hline
\end{tabular}
\end{center}

\nl
The first equation gives
\begin{eqnarray}\label{commTH}
T_{j,\b}h_{\b,i} &=& 
h_{\a_j,i} T_{j,\b}
\end{eqnarray}
and the second
\begin{eqnarray}\label{commlaweq}
T_{j,\b} &=& h_{\a_j,i}^{-1}T_{j,\b-\a_i}.
\end{eqnarray}

The third case gives an equation that is equivalent to (\ref{commlaweq}).
The fourth equation is part of (\ref{knotleq}) in Table \ref{algtable}
(namely the part where $j\not \sim i$).

\np
Next, we assume $i\sim j$.
A practical rule is
\begin{eqnarray*}\label{alp2alpteq}
\tau_i\tau_jx_{\a_i}
&=&
\tau_i(-mx_{\a_i} + x_{\a_i+\a_j}) = 
x_{\a_j}.
\end{eqnarray*}
We distinguish cases according to the values of
$(\a_i,\b)$ and $(\a_j,\b)$. Since each inner product, for distinct 
roots
is one of $1$, $0$, $-1$, there are six cases to consider up to 
interchanges of $i$ and $j$.
However, as in the proof of Lemma \ref{monmorphismlm} for $i\sim j$, the case 
$(\a_i,\b)=(\a_j,\b)= -1$ does not occur.

For the sake of brevity, let us denote the images of the left hand side 
and
the right hand side of (\ref{genTrepeqadj}) 
on $x_\b$ by
$LHS$ and $RHS$, respectively. 

\nl
Case $(\a_i,\b)=(\a_j,\b)= 1$. Then
$(r_i\b,\a_j)=(\b-\a_i,\a_j)=2$, so
$\b = \a_i+\a_j$.
Now
\begin{eqnarray*}
RHS& =&
x_{\a_j}(T_{i,\b} -mT_{j,\a_j}) + x_{\b}T_{j,\a_j}.
\end{eqnarray*}
Comparison with the same expression but then
$j$ and $i$ interchanged yields $LHS$. This leads to the 
equations
$T_{i,\b} = m T_{j,\a_j}$ and $T_{j,\a_j} = T_{i,\a_i}$.
In view of the latter, and connectedness of the diagram
there is an element $z$ in $Z_0^{(0)}$ such that
\begin{eqnarray}\label{11eq2}
T_{i,\a_i} &=& z \ \ \mbox{ for all } \ \ i.
\end{eqnarray}
Consequently, the former equation reads
\begin{eqnarray}\label{11eq1}
T_{i,\b} &=& m z.
\end{eqnarray}
By the requirement
$T_{i,\a_i} =1$ in the hypotheses, we must have $z = 1$.

\nl
Case $(\a_i,\b)=(\a_j,\b)= 0$. This gives
\begin{eqnarray*}
RHS& =&
x_{\a_j}(T_{i,\b} -m T_{j,\b}h_{\b,i} ) +
x_{\a_j+\a_i} T_{j,\b} h_{\b,i} + x_{\a_i} T_{i,\b} 
h_{\b,j}h_{\b,i}
\end{eqnarray*}
and
$LHS$ can be obtained from the above
by interchanging the indices $i$ and $j$.
Comparison of each of the coefficients of
$x_{\a_i}$, $x_{\a_j+\a_i}$, $x_{\a_j}$ gives
\begin{eqnarray}\label{00eq}
T_{i,\b} h_{\b,j} & = & T_{j,\b} h_{\b,i} 
\ \ \mbox{if} \ \ (\a_i,\b)=(\a_j,\b)= 0
\ \ \mbox{and} \ \ (\a_i,\a_j) = -1.
\end{eqnarray}

\np
Since the other cases come down to similar computations, we only list
the results.

\bigskip\noindent
Case $(\a_i,\b)=0$, $(\a_j,\b)= -1$.
Here we have
\begin{eqnarray*}
RHS&=&
x_{\a_i} (-m T_{i,\b} h_{\b,i} + T_{i,\b+\a_j} h_{\b,i})
+ x_{\a_{j}} (-m T_{j,\b} h_{\b,i} + T_{i,\b})\\&&
+ x_{\a_i+\a_j} (T_{j,\b} h_{\b,i})
\end{eqnarray*}
and
\begin{eqnarray*}
LHS
&=&
x_{\a_i} ( m^2 T_{i,\b} + T_{j,\b} -m T_{i,\b+\a_j}) 
+ x_{\a_j} (-m T_{j,\b}h_{\b,i} -m T_{j,\b+\a_j} \\&&
+T_{j,\b+\a_j+\a_i}) 
+ x_{\a_i+\a_j} (-m T_{i,\b} + T_{i,\b+\a_j}),\\
\end{eqnarray*}
which gives

\begin{eqnarray}\label{0-1eq1}
T_{i,\b+\a_j} &=& T_{j,\b} h_{\b,i} + m T_{i,\b}, \\
T_{j,\b+\a_j+\a_i} &=& T_{i,\b} + mT_{j,\b+\a_j}.\label{0-1eq2}
\end{eqnarray}

\bigskip\noindent
Case $(\a_i,\b)=0$, $(\a_j,\b)= 1$.
\begin{eqnarray*}
RHS
&=&
x_{\a_i} T_{i,\b-\a_j} h_{\b,i} +
x_{\a_j} (-mT_{j,\b} h_{\b,i} + T_{i,\b}) +
x_{\a_j+\a_i} T_{j,\b} h_{\b,i}
\end{eqnarray*}
and
\begin{eqnarray*}
LHS
&=&
x_{\a_i} (T_{j,\b} - mT_{i,\b-\a_j}) + 
x_{\a_j} T_{j,\b-\a_j-\a_i} +
x_{\a_i+\a_j} T_{i,\b-\a_j}
\end{eqnarray*}
whence
\begin{eqnarray}\label{01eq1}
T_{i,\b} &=& T_{j,\b-\a_j-\a_i} + mT_{i,\b-\a_j},
\\
T_{j,\b} &=& T_{i,\b-\a_j} h_{\b,i}^{-1}.\label{01eq2}
\end{eqnarray}

\nl
Case $(\a_i,\b)=1$, $(\a_j,\b)= -1$.
Now
\begin{eqnarray*}
RHS 
&=&
x_{\a_i} (T_{i,\b-\a_i} h_{\b-\a_i,j}) +
x_{\a_j} (T_{i,\b} - m T_{j,\b-\a_i}) +
x_{\a_i+\a_j} (T_{j,\b-\a_i}) \\
\end{eqnarray*}
and
\begin{eqnarray*}
LHS&=&
x_{\a_i} (m^2 T_{i,\b} - m T_{i,\b+\a_j} + T_{j,\b}) + 
x_{\a_j} (T_{j,\b+\a_j} h_{\b+\a_j,i} - m T_{j,\b-\a_i}) +\\
&& x_{\a_i+\a_j} (T_{i,\b+\a_j} - m T_{i,\b}) \\
\end{eqnarray*}
whence
%
\begin{eqnarray}\label{1-1eq1}
T_{j,\b} &=& T_{i,\b-\a_i} h_{\b-\a_i,j} + m T_{j,\b-\a_i}, \\
T_{j,\b+\a_j} &=& T_{i,\b}h_{\b+\a_j,i}^{-1},\label{1-1eq2} \\
T_{i,\b+\a_j} &=& T_{j,\b-\a_i} + mT_{i,\b}.
\label{1-1eq3}
\end{eqnarray}

\begin{center}
\begin{table}
\caption{Equations for $T_{i,\b}$}
\label{algtable}
\begin{tabular}{|l|c|c|}
\hline
$T_{i,\b}$&condition&reference\\
\hline
$0$& $\b=\a_j$ and $i\ne j$&(\ref{knotleq}) \\
\hline
$1$&$\b=\a_i$&(\ref{11eq2}) \\
\hline
$m$&$\b=\a_i+\a_j$&(\ref{11eq1}) \\
\hline
$h_{\a_i,j}^{-1}T_{i,\b-\a_j}$
&$(\a_j,\b)=1$ and $(\a_i,\a_j)=0$
&(\ref{commlaweq})\\
\hline
$T_{j,\b-\a_i-\a_j}+mT_{i,\b-\a_j}
$&$(\a_i,\b)=0$ and $(\a_j,\b)=1$&(\ref{01eq1})
\\
& and $(\a_i,\a_j)=-1$&\\
\hline
$T_{j,\b-\a_j}h_{\b-\a_j,i}+ m T_{i,\b-\a_j}
$&$(\a_i,\b)=-1 $ and $(\a_j,\b)=1$& (\ref{1-1eq1}) 
\\
& and $(\a_i,\a_j)=-1$&\\
\hline
$T_{j,\b-\a_i}h_{\b,j}^{-1}$&$(\a_i,\b)=1$ and 
$(\a_j,\b)=0$&(\ref{01eq2}) 
\\
& and $(\a_i,\a_j)=-1$&\\
\hline
\end{tabular}
\end{table}
\end{center}

\np
We now consider the coefficients of $l^{-2}$ and of $l^{-3}$
in the equations (B1), (B2) for $\s_i$.
We claim that, given (\ref{commlaweq})--(\ref{1-1eq3}),
a necessary condition for the corresponding
equations to hold is
\begin{eqnarray}\label{knotleq}
T_{k,\a_j} &=& 0 \quad\ \mbox{if }\ k\ne j.
\end{eqnarray}
To see this, note that,
if $k\not\sim j$,
the coefficient of $l^{-2}$ gives
$T_kT_j=T_jT_k$
which, applied to $x_{\a_j}$,
yields (\ref{knotleq}). If $k\sim j$, note 
$$
T_k\tau_jx_{\a_k} = 
T_k(-mx_{\a_k}+x_{\a_k+\a_j})
=0$$ as $T_{k,\a_k+\a_j}=mz = mT_{k,\a_k}$ by (\ref{11eq1}). Now use the action of
$$T_j\tau_k T_j + \tau_jT_kT_j+T_jT_k\tau_j= T_k\tau_j T_k + 
\tau_k
T_jT_k+T_kT_j\tau_k 
$$
on $x_{\a_j}$. 
We see only the middle terms do not vanish because of the relation 
above
and so 
$$
\tau_jx_{\a_k}T_{k,\a_j}z=\tau_kx_{\a_j} T_{j,\a_k}T_{k,\a_j}.
$$
By
considering the coefficient of $x_{\a_k}$, which occurs only on the
left hand side, we see that (\ref{knotleq}) holds.

A consequence of this is that $T_iT_j=0$ if $i\neq j$. Now all 
the
equations for the $l^{-2}$ and $l^{-3}$ coefficients are easily satisfied.
In the noncommuting case of $l^{-2}$, the first
terms on either side are $0$ by the relation
above and the other terms are $0$ as $T_jT_k=0$.

\np
We have seen that, in order for $s_i\mapsto \sigma_i$ to determine a representation,
the $T_{i,\b}$ have to satisfy the equations 
(\ref{commTH})--(\ref{knotleq}).
This system of equations, however,
is redundant. Indeed, when the root in the index of
the left hand side of 
(\ref{0-1eq1}) is set to $\c$, we obtain 
(\ref{1-1eq1}) for $\c$ instead of $\b$.
Similarly, (\ref{0-1eq2}) is equivalent to
(\ref{01eq1}), while
(\ref{1-1eq2}) is equivalent to
(\ref{01eq2}),
and (\ref{1-1eq3}) is equivalent to
(\ref{01eq1}).
Consequently,
in order to finish the proof that Table \ref{algtable} contains a sufficient
set of relations, we must show that (\ref{00eq}) and (\ref{commTH}) follow
from those of the table. These proofs are given in Lemmas \ref{00eq-lm} and
\ref{Conj-rel-Lm} below.

It remains to establish that the matrices $\s_k$ are invertible. To prove
this, we observe that the linear transformation $\s_k^2+m\s_k-\ident$ maps $V$
onto the submodule spanned by $x_{\a_k}$ and that the image of $x_{\a_k}$
under $\s_k$ is $x_{\a_k}l^{-1}$.  This is easy to establish and will be shown 
in Lemma \ref{eiproj} below.
\end{proof}

\begin{Lm}\label{00eq-lm}
The equations in (\ref{00eq}) are consequences of the relations of Table \ref{algtable}.  
\end{Lm}

\begin{proof}
The equation says that
$T_{k,\b}h_{\b,j} = T_{j,\b}h_{\b,k}$ whenever $(\a_k,\b)=(\a_j,\b)=0$ and
$k\sim j$. We prove this by induction on the height of $\b$. The initial case
of $\b$ having height~$1$ is direct from (\ref{knotleq}).  Suppose therefore,
$\het(\b)>1$.  There exists $m\in\{1,\ldots,n\}$ such that $(\a_m,\b) = 1$.
If $(\a_m,\a_k) = (\a_m,\a_j)=0$, then, by the induction hypothesis and
(\ref{H-orth-rel}), $T_{k,\b-\a_m}h_{\b,j} = T_{k,\b-\a_m}h_{\b-\a_m,j} =
T_{j,\b-\a_m}h_{\b-\a_m,k} T_{j,\b-\a_m}h_{\b,k}$, so, applying
(\ref{commlaweq}) twice, we find 
$$T_{k,\b}h_{\b,j} =
h_{\a_k,m}^{-1}T_{k,\b-\a_m}h_{\b,j} = h_{\a_j,m}^{-1}T_{j,\b-\a_m} h_{\b,k}=
T_{j,\b} h_{\b,k},$$ 
as required.

Therefore, interchanging $k$ and $j$ if necessary, we may assume 
that
$j\sim m$, whence $k\not\sim m$ (as the Dynkin diagram
contains no triangles).
Now $\delta=\b-\a_m-\a_j$ and $\c=\delta-\a_k$ are positive roots and
$(\a_k,\delta)=1$, so (\ref{commlaweq}) gives
$T_{m,\c} =  h_{\a_m,k}T_{m,\delta}$, which, by induction on height,
and (\ref{H-dist2-rel}),
leads to
$ h_{\a_k,m}^{-1}T_{j,\c}  =h_{\a_m,k}^{-1} T_{m,\c} h_{\c,j}h_{\c,m}^{-1} = T_{m,\delta}  h_{\c,j}h_{\c,m}^{-1}$.
Observing that, by straightforward application of the braid relations and the
definition of $h_{\b,k}$, we also have
\begin{eqnarray*}
 h_{\c,j}h_{\c,m}^{-1}h_{\b,j}&=& h_{\b,k}\\
h_{\delta,m}^{-1}h_{\b,j} &=& h_{\b-\a_m,k}^{-1}h_{\b,k}
\end{eqnarray*}
we derive
\begin{eqnarray*}
T_{k,\b}h_{\b,j} &=&  h_{\a_k,m}^{-1}T_{k,\b-\a_m} h_{\b,j} \qquad 
\hfill\phantom{mhhhhhhhhi}\,
\mbox{ by (\ref{commlaweq})}\\
&=& h_{\a_k,m}^{-1}
\left(T_{j,\c} + m T_{k,\delta}\right) h_{\b,j}
\phantom{mxhhhhhh}\, \mbox{ by (\ref{01eq1})}\\
&=& T_{m,\delta}  h_{\c,j}h_{\c,m}^{-1}h_{\b,j} + m  T_{k,\delta}h_{\delta,m}^{-1} h_{\b,j}\qquad \hfill
\phantom{much} \mbox{ by the above and (\ref{commTH})}\\
&=& T_{m,\delta}T_{m,\delta}  h_{\b,k} + m  T_{k,\delta}h_{\b-\a_m,k}^{-1} h_{\b,k}\qquad \hfill
\phantom{much} \mbox{ by the above }\\
&=& \left(T_{m,\delta} +mT_{j,\b-\a_m}\right)h_{\b,k} \qquad \hfill
\phantom{much} \mbox{ by (\ref{01eq2})}\\
&=& T_{j,\b}h_{\b,k}, \qquad \hfill \phantom{mxucccchmorrejjjjj}
\mbox{ by (\ref{01eq1})}
\end{eqnarray*}
as required.
\end{proof}

\np
The relation (\ref{commTH}) is new compared to \cite{CW}.
But it is superfluous. In order to see this, we first prove some auxiliary
claims.

\begin{Lm}\label{H-extra-rel-Lm}
Let $h$, $k$ be generators (or conjugates thereof) in the Hecke algebra $Z_0^{(0)}$.
Then, for any $t\in Z_0^{(0)}$, 
\begin{enumerate}
\item $h^{-1}t - tk^{-1} = ht - tk$,
\item $h^{-1} (t+h^{-1}tk^{-1})k = t + h^{-1} t k^{-1}$.
\end{enumerate}
\end{Lm}

\begin{proof}
(i). Expand the left hand side and use that $z^{-1} = z + m$ for every
conjugate of a generator.

\nl(ii). By (i),
$tk+h^{-1}t = ht + tk^{-1}$. Multiplying both sides from the left by $h^{-1}$
and pulling out a factor $k$ at the right of the left hand side, we find the required relation.
\end{proof}

\begin{Lm}\label{Conj-rel-Lm}
The relation (\ref{commTH}) is a consequence of those of Table \ref{algtable}.
\end{Lm}

\begin{proof}
Suppose that the positive root $\b$ and the distinct nodes $l$, $i$ satisfy 
$(\a_l,\b) = 0$ and $i\not\sim l$. We want to establish
(\ref{commTH}) by induction on $\het(\b)$.
By Corollary \ref{jastline}(i), we know that $T_{i,\b} = 0$ if
$i\not\in\Supp(\b)$, so we need only consider cases where $i\in\Supp(\b)$.

If $\het(\b) = 1$, then,
by (\ref{11eq2}) and (\ref{knotleq}), $T_{i,\b} = 0$  and there is nothing to prove unless $\b
= \a_i$. In the latter case $T_{i,\b} = 1$ and
$h_{\a_i,l}^{-1}T_{i,\b}h_{\b,l} = h_{\a_i,l}^{-1}h_{\a_i,l} = 1$, so
(\ref{commTH}) is satisfied.

If $\het(\b) = 2$, then $\b= \a_i+\a_j$ for some $j$ and $T_{i,\b} = m$ by (\ref{11eq1}).
As $\a_l$ is orthogonal to both $\b$ and $\a_i$, it must be orthogonal to
$\a_j$ as well. 
Now $h_{\a_i,l}^{-1}T_{i,\b}h_{\b,l} = m h_{\a_i,l}^{-1}h_{\a_i+\a_j,l} = 
m d_{\a_i}^{-1}s_l^{-1}d_{\a_i}d_{\a_i}^{-1}s_j^{-1}s_ls_jd_{\a_i}
= m$, as required.

\np
Case (\ref{commlaweq}): there is a node $j$ with
$(\a_j,\b)=1$ and $(\a_i,\a_j)=0$. Then
$T_{i,\b} = h_{\a_i,j}^{-1}T_{i,\b-\a_j}$.
If $j\not\sim l$, we find
\begin{eqnarray*}
h_{\a_i,l}^{-1}T_{i,\b}h_{\b,l} &= &
h_{\a_i,l}^{-1}
h_{\a_i,j}^{-1}T_{i,\b-\a_j}
h_{\b,l} \qquad \hfill \phantom{mxucchmucccchmorre}
\mbox{ by (\ref{commlaweq})} \\
&= &
h_{\a_i,j}^{-1}h_{\a_i,l}^{-1}
T_{i,\b-\a_j}
h_{\b-\a_j,l} \qquad \hfill \phantom{mucccchmorre}
\mbox{ by (\ref{H-com-rel}) and (\ref{H-orth-rel})}\\
&= &
h_{\a_i,j}^{-1}T_{i,\b-\a_j}
\qquad \hfill \phantom{mxucchmucccchmorre}
\mbox{ by induction} \\
&= &
T_{i,\b}. \qquad \hfill \phantom{mxucchmucccchmorre}
\mbox{ by (\ref{commlaweq})} 
\end{eqnarray*}

If $j\sim l$, we find
\begin{eqnarray*}
h_{\a_i,l}^{-1}T_{i,\b}h_{\b,l} &= &
h_{\a_i,l}^{-1}
h_{\a_i,j}^{-1}T_{i,\b-\a_j}
h_{\b,l}  \qquad \hfill \phantom{mxucchmucccchmorre}
\mbox{ by (\ref{commlaweq})} 
\\
&= &
h_{\a_i,l}^{-1}h_{\a_i,j}^{-1}h_{\a_i,l}^{-1}
T_{i,\b-\a_j-\a_l}
h_{\b,l} \qquad \hfill \phantom{mxucchmucccchmorre}
\mbox{ by (\ref{commlaweq})} 
 \\
&= &
h_{\a_i,j}^{-1}h_{\a_i,l}^{-1}h_{\a_i,j}^{-1}
T_{i,\b-\a_j-\a_l}
h_{\b-\a_j-\a_l,j}  \qquad \hfill \phantom{e}
\mbox{ by (\ref{H-braid-rel}) and (\ref{H-10-rel}})
\\
&= &
h_{\a_i,j}^{-1}h_{\a_i,l}^{-1}
T_{i,\b-\a_j-\a_l}
 \qquad \hfill \phantom{mxucchmucccchmorre}
\mbox{ by induction} 
 \\
&= &
T_{i,\b}.
 \qquad \hfill \phantom{mxucchmucccchmorre}
\mbox{ by (\ref{commlaweq}) applied twice} 
\end{eqnarray*}
This ends Case (\ref{commlaweq}).

\np
Case (\ref{01eq1}): $(\a_i,\b)=0$ and 
there is a node $j\sim i$ with $(\a_j,\b)=1$.
Then $T_{i,\b} = T_{j,\b-\a_i-\a_j}+mT_{i,\b-\a_j}$.
Now $h_{\a_i,l}^{-1}T_{i,\b}h_{\b,l} =
h_{\a_i,l}^{-1}(T_{j,\b-\a_i-\a_j}+mT_{i,\b-\a_j})h_{\b,l}$.

If $j\not\sim l$, we find
\begin{eqnarray*}
h_{\a_i,l}^{-1}T_{i,\b}h_{\b,l} &= &
h_{\a_i,l}^{-1}(T_{j,\b-\a_i-\a_j}+mT_{i,\b-\a_j})h_{\b,l} 
 \hfill \phantom{mucccchmorre}
\mbox{ by (\ref{01eq1})} \\
&= &
h_{\a_i,l}^{-1}T_{j,\b-\a_i-\a_j}h_{\b-\a_i-\a_j,l} +mh_{\a_i,l}^{-1}T_{i,\b-\a_j}h_{\b-\a_j,l} 
 \quad \hfill \phantom{}
\mbox{ by (\ref{H-orth-rel})}\\
&= &
T_{j,\b-\a_i-\a_j} +mT_{i,\b-\a_j}
\qquad \hfill \phantom{mxucchmucccchmorre}
\mbox{ by induction} \\
&= &
T_{i,\b}. \qquad \hfill \phantom{mxucchmucccchmorre}
\mbox{ by (\ref{01eq1})} 
\end{eqnarray*}

If $j\sim l$, we claim
\begin{eqnarray}\label{delta27-eq}
T_{i,\b}&=&
T_{l,\delta} + m (T_{j,\c} + h_{\a_i,l}^{-1}T_{j,\c}h_{\b,l}^{-1}),
\end{eqnarray}
where $\c=\b-\a_i-\a_j-\a_l$ and where $\delta=\c-\a_j$ are positive roots.
For
\begin{eqnarray*}
T_{i,\b}&=&
T_{j,\b-\a_i-\a_j}+mT_{i,\b-\a_j}
 \hfill \phantom{mucccchmorre}
\mbox{ by (\ref{01eq1})} \\
&= &
(T_{l,\delta}+mT_{j,\c})
+ m 
h_{\a_i,l}^{-1}T_{i,\b-\a_j-\a_l}
 \hfill \phantom{mucccchmorre}
\mbox{ by (\ref{01eq1}) and (\ref{commlaweq})} \\
&= &
T_{l,\delta}+mT_{j,\c}
+ m 
h_{\a_i,l}^{-1}T_{j,\c}h_{\b-\a_j-\a_l,j}^{-1}
 \hfill \phantom{mucccchmorre}
\mbox{ by (\ref{01eq2})} \\
&= &
T_{l,\delta}+mT_{j,\c}
+ m h_{\a_i,l}^{-1}T_{j,\c}h_{\b,l}^{-1}.
 \phantom{mucccchmorre}
\mbox{ by (\ref{H-10-rel})} 
\end{eqnarray*}
By (\ref{H-10-rel}), we have
$h_{\b,l} = h_{\b-\a_j-\a_l,j} = h_{\delta,i}$, so, by induction we find
$h_{\a_i,l}^{-1}T_{l,\delta} h_{\b,l} = 
T_{l,\delta} h_{\delta,i}^{-1})h_{\b,l} = T_{l,\delta} $.
So the first summand of (\ref{delta27-eq}) is
invariant under under simultaneous
left multiplication by
$h_{\a_i,l}^{-1}$ and right multiplication by  $h_{\b,l}$.
The same holds for the second summand, 
$m(T_{j,\c} + h_{\a_i,l}^{-1}T_{j,\c}h_{\b,l}^{-1})$
by  Lemma \ref{H-extra-rel-Lm} applied with
$h = h_{\a_i,l}$,  $k = h_{\b,l}$, and $t = T_{j,\c}$.
Consequently (\ref{commTH}) holds for $T_{i,\b}$ in Case (\ref{01eq1}).

\np
Case (\ref{1-1eq1}): $(\a_i,\b)=-1$ and 
there is a node $j\sim i$ with $(\a_j,\b)=1$.
Then
$T_{i,\b} = T_{j,\b-\a_j}h_{\b-\a_j,i}+ m T_{i,\b-\a_j}$.
Now $h_{\a_i,l}^{-1}T_{i,\b}h_{\b,l} =
h_{\a_i,l}^{-1}(T_{j,\b-\a_j}h_{\b-\a_j,i}+ m T_{i,\b-\a_j})h_{\b,l}$.

If $j\not\sim l$, we find
\begin{eqnarray*}
h_{\a_i,l}^{-1}T_{i,\b}h_{\b,l} &= &
h_{\a_i,l}^{-1}(T_{j,\b-\a_j}h_{\b-\a_j,i}+ m T_{i,\b-\a_j})h_{\b,l}
\hfill \phantom{mhmucccchmorre}
\mbox{ by (\ref{01eq1})} \\
&= &
h_{\a_i,l}^{-1}T_{j,\b-\a_j}h_{\b-\a_j,l}h_{\b-\a_j,i} +mh_{\a_i,l}^{-1}T_{i,\b-\a_j}h_{\b-\a_j,l} 
 \quad \hfill \phantom{}
\mbox{ by (\ref{H-orth-rel}) and (\ref{H-com-rel})}\\
&= &
T_{j,\b-\a_j}h_{\b-\a_j,i} +mT_{i,\b-\a_j}
\qquad \hfill \phantom{mxucchmucccchmorre}
\mbox{ by induction} \\
&= &
T_{i,\b}. \qquad \hfill \phantom{mxucchmucccchmorre}
\mbox{ by (\ref{01eq1})} 
\end{eqnarray*}

If $j\sim l$, we claim
\begin{eqnarray}\label{gamma36-eq}
T_{i,\b}&=&
T_{l,\c}h_{\c,j}h_{\b-\a_j,i} + m\left(T_{j,\c}h_{\b-\a_j,i} + h_{\a_i,l}^{-1}T_{i,\c}\right),
\end{eqnarray}
where $\c=\b-\a_j-\a_l$ is a positive root.
For
\begin{eqnarray*}
T_{i,\b}&=&
T_{j,\b-\a_j}h_{\b-\a_j,i}+mT_{i,\b-\a_j}
 \hfill \phantom{mucccchmorre}
\mbox{ by (\ref{1-1eq1})} \\
&= &
T_{l,\c}h_{\c,j}h_{\b-\a_j,i}
+mT_{j,\c}h_{\b-\a_j,i}
+mh_{\a_i,l}^{-1}T_{i,\c}.
 \phantom{mucccchmorre}
\mbox{ by (\ref{1-1eq1}) and (\ref{commlaweq})}
\end{eqnarray*}
By Lemma \ref{H-rel}, we 
have
\begin{eqnarray*}
&&h_{\c,i}^{-1}h_{\c,j}h_{\b-\a_j,i}h_{\b,l} =
d_{\b}^{-1}(s_j^{-1}s_l^{-1}s_i^{-1}s_ls_j)
(s_j^{-1}s_l^{-1}s_js_ls_j)
(s_j^{-1}s_is_j)
(s_l)d_{\b}\\
&&= d_{\b}^{-1}
(s_j^{-1}s_i^{-1}s_j)
(s_j^{-1}s_l^{-1}s_js_ls_j)
(s_is_js_i^{-1})
(s_l)
d_{\b}
= 
d_{\b}^{-1}
s_j^{-1}s_i^{-1}s_l^{-1}s_js_ls_is_j
s_l
d_{\b}
\\
&&= 
d_{\b}^{-1}
s_j^{-1}s_l^{-1}s_i^{-1}s_js_is_ls_j
s_l
d_{\b}
= 
d_{\b}^{-1}
s_j^{-1}s_l^{-1}s_js_is_j^{-1}s_ls_j
s_l
d_{\b}\\
&& = 
d_{\b}^{-1}
s_ls_j^{-1}s_l^{-1}s_is_ls_j
d_{\b}
= 
d_{\b}^{-1}
s_ls_j^{-1}s_is_j
d_{\b}
\\&& = h_{\c,j}h_{\b-\a_j}.
\end{eqnarray*}
Hence, using induction, we find for the first summand of (\ref{gamma36-eq})
$$h_{\a_i,l}^{-1}\left(T_{l,\c} h_{\c,j}h_{\b-\a_j,i}\right) h_{\b,l} = 
T_{l,\c} h_{\c,i}^{-1}h_{\c,j}h_{\b-\a_j,i}
h_{\b,l} = T_{l,\delta}h_{\c,j}h_{\b-\a_j},$$
proving that it is
invariant under under simultaneous
left multiplication by
$h_{\a_i,l}^{-1}$ and right multiplication by  $h_{\b,l}$.

The same holds for the second summand, 
$m(T_{j,\c}h_{\b-\a_j,i} + h_{\a_i,l}^{-1}T_{i,\c})$
as we shall establish next.
First of all, note that $h_{\c,j} = h_{\b,l}$ by (\ref{H-10-rel})
and that $h_{\c,i} = h_{\b-\a_j,i}$ by (\ref{H-orth-rel}). 
Moreover, by Corollary \ref{jastline}(iii),
$T_{i,\c}h_{\c,j} = T_{j,\c}h_{\c,i} $.
Substituting all this in 
the second summand, we obtain
\begin{eqnarray*}
m(T_{j,\c}h_{\b-\a_j,i} + h_{\a_i,l}^{-1}T_{i,\c}) & = &
m(T_{j,\c}h_{\c,i} + h_{\a_i,l}^{-1}T_{i,\c})
=m(T_{i,\c}h_{\c,j} + h_{\a_i,l}^{-1}T_{i,\c})
\\
&=&m(T_{i,\c}h_{\b,l} + h_{\a_i,l}^{-1}T_{i,\c}).
\end{eqnarray*}
Again, using Lemma \ref{H-extra-rel-Lm} applied with
$h = h_{\a_i,l}$,  $k = h_{\b,l}$, and $t = T_{i,\c}$,
we find the required invariance.
Consequently (\ref{commTH}) holds for $T_{i,\b}$ in Case (\ref{01eq1}).

\np
Case (\ref{01eq2}): $(\a_i,\b)=1$ and 
there is a node $j\sim i$ with $(\a_j,\b)=0$.
Then
$T_{i,\b} = T_{j,\b-\a_i}h_{\b,j}^{-1}$.
Now $h_{\a_i,l}^{-1}T_{i,\b}h_{\b,l} =
h_{\a_i,l}^{-1}T_{j,\b-\a_i}h_{\b,j}^{-1}h_{\b,l}$.

If $j\not\sim l$, we find
\begin{eqnarray*}
h_{\a_i,l}^{-1}T_{i,\b}h_{\b,l} &= &
h_{\a_i,l}^{-1}T_{j,\b-\a_i}h_{\b,j}^{-1}h_{\b,l} 
\qquad \hfill \phantom{mxucchmucccchmorre}
\mbox{ by (\ref{01eq2})} \\
&= &
h_{\a_i,l}^{-1}T_{j,\b-\a_i}h_{\b-\a_i,l}h_{\b,j}^{-1} 
 \quad \hfill \phantom{}
\mbox{ by (\ref{H-com-rel}) and (\ref{H-orth-rel})}\\
&= &
T_{j,\b-\a_i}h_{\b,j}^{-1} 
\qquad \hfill \phantom{mxucchmucccchmorre}
\mbox{ by induction} \\
&= &
T_{i,\b}. \qquad \hfill \phantom{mxucchmucccchmorre}
\mbox{ by (\ref{01eq2})} 
\end{eqnarray*}

If $j\sim l$, observe that
$h_{\b-\a_i,l}^{-1}h_{\b,j}^{-1}h_{\b,l}
= h_{\b-\a_i-\a_j,i}  h_{\b-\a_i,l} ^{-1}    h_{\b,j}^{-1}$ in view of (\ref{H-orth-rel}), (\ref{H-10-rel}), and
(\ref{H-braid-rel}). 
Also, $h_{\a_i,l} = h_{\a_l,i}$ by a double application of (\ref{H-0-1-rel}).
Therefore,
\begin{eqnarray*}
h_{\a_i,l}^{-1}T_{i,\b}h_{\b,l} &= &
h_{\a_i,l}^{-1}T_{l,\b-\a_i-\a_j}h_{\b-\a_i,l}^{-1}h_{\b,j}^{-1}h_{\b,l} 
 \quad \hfill \phantom{}
\mbox{ by (\ref{01eq2}) twice}\\
&= &
h_{\a_l,i}^{-1}T_{l,\b-\a_i-\a_j} h_{\b-\a_i-\a_j,i}  h_{\b-\a_i,l} ^{-1}h_{\b,j}^{-1}
 \phantom{more}
\mbox{ by the above} \\
&= &
T_{l,\b-\a_i-\a_j}  h_{\b-\a_i,l} ^{-1}h_{\b,j}^{-1}
\qquad \hfill \phantom{mxucchmucccchmorre}
\mbox{ by induction} \\
&= &
T_{i,\b}. \qquad \hfill \phantom{mxucchmucccchmorre}
\mbox{ by (\ref{01eq2}) twice} 
\end{eqnarray*}
\end{proof}

\begin{Cor}\label{jastline}
If the $T_{i,\b}\in Z_0^{(0)}$ satisfy the equations in Table
\ref{algtable}, then these obey the following rules,
where $\het(\b)$ stands for the height of $\b$.
\begin{enumerate}[(i)]
\item $T_{i,\b}=0$ whenever $i\not\in \Supp(\b)$.
\item If $(\a_i,\b) = 1$, then
$T_{i,\b} = m d_{\a_i}^{-1}s_\b^{-1}s_is_\b d_{\b}$.
\item 
If $(\a_i,\b)=(\a_j,\b)=0$
and $(\a_i,\a_j)=-1$, then $T_{i,\b} h_{\b,j}= T_{j,\b} h_{\b,i}$.
\end{enumerate}
\end{Cor}

\np
\begin{proof}
(i) follows from (\ref{knotleq}) by use of (\ref{commlaweq}) and
(\ref{1-1eq1}). Observe that, if $i\not\in\Supp(\b)$ and $(\a_j,\b) 
=1$
for some $j\sim i$, then $j\not\in\Supp(\b-\a_j)$.

\nl(ii).  By induction on $\het(\b)$. The assertion is vacuous when $\het(\b)
 =1$.
Suppose $\het(\b) = 2$. Then $s_\b = s_js_is_j$ for some node $j$ adjacent to
 $i$ in $M$.
Therefore, $m d_{\a_i}^{-1}s_\b^{-1}s_is_\b d_{\b} = m
 d_{\a_i}^{-1}s_j^{-1}s_i^{-1}s_j^{-1}s_is_js_is_j d_{\b}
= m d_{\b}^{-1}s_j^{-1}s_j^{-1}s_i^{-1}s_j^{-1}s_is_js_is_j d_{\b}
= m$ and, by (\ref{11eq1}) $ T_{i\,\b} = m$, as required.

Now suppose $\het(\b)>2$.

If $j$ is a node distinct from $i$ such that $(\a_j,\b) = 1$, then, necessarily, $i\not\sim j$
(for otherwise $(\a_i,\b-\a_j) = 2$, so $\b = \a_i+\a_j$, contradicting
$\het(\b) > 2$).
Now
(\ref{commlaweq}) applies, giving
\begin{eqnarray*}
T_{i,\b} 
 &=&
h_{\a_i,j}^{-1} T_{i,\b-\a_i}
\qquad \hfill \phantom{mxucchmucccchmorre}
\mbox{by (\ref{commlaweq})} \\
 &=&
md_{\a_i}^{-1}s_j^{-1}s_{\b-\a_j}^{-1}s_is_{\b-\a_j}s_jd_{\b}
\qquad \hfill \phantom{mxucchmucccchmorre}
\mbox{by induction} \\
 &=&
md_{\a_i}^{-1} s_{\b}^{-1}   s_js_is_j^{-1}s_{\b}d_{\b}
\qquad \hfill \phantom{mxucchmucccchmorre}
\mbox{by definition of } s_\b \\
 &=&
md_{\a_i}^{-1} s_{\b}^{-1}   s_is_{\b}d_{\b},
\qquad \hfill \phantom{mxucchmucccchmorre}
\mbox{as } s_is_j = s_js_i 
\end{eqnarray*}
as required.

Suppose $l$ is a node distinct from $i$ such that $(\a_l,\b) = 0$ and $i\sim l$.
Then
(\ref{01eq2}) applies, giving
\begin{eqnarray*}
T_{i,\b} 
 &=&
T_{l,\b-\a_i}h_{\b,l}^{-1}
\qquad \hfill \phantom{mxucchmucgggccchmorre}
\mbox{by (\ref{01eq2})} \\
 &=&
m d_{\a_l}^{-1}(s_{\b-\a_i}^{-1}s_l)s_{\b-\a_i}(d_{\b-\a_i}d_{\b}^{-1})s_l^{-1}d_{\b}
\qquad \hfill \phantom{mxucch}
\mbox{by induction} \\
 &=&
m d_{\a_i}^{-1}(s_ls_i^{-1}s_l^{-1})
s_{\c}^{-1} s_ls_{\c}
(s_ls_is_l^{-1})
d_{\b}
\qquad \hfill \phantom{mc}
\mbox{by definition of } d_\b \mbox{ and } s_\b \\
 &=&
m d_{\a_i}^{-1}s_i^{-1}
s_l^{-1}s_is_{\c}^{-1} 
s_ls_{\c}s_i^{-1}s_l
s_i
d_{\b}
\qquad \hfill \phantom{mhmuc}
\mbox{by the braid relation} 
\\
 &=&
m d_{\a_i}^{-1}s_i^{-1}
s_l^{-1}s_{\c}^{-1} s_i
s_ls_i^{-1}s_{\c}s_l
s_i
d_{\b}
\qquad \hfill \phantom{mcchmuch}
\mbox{by Lemma \ref{isb=sbi-lm}} 
\\
 &=&
m d_{\a_i}^{-1} (s_i^{-1} s_l^{-1}s_{\c}^{-1} s_l^{-1})
(s_is_ls_{\c}s_ls_i)
d_{\b}
\qquad \hfill \phantom{mh}
\mbox{by the braid relation} 
\\
 &=&
m d_{\a_i}^{-1}s_{\b}^{-1} s_i s_{\b} d_{\b}
\qquad \hfill \phantom{mmuch}
\mbox{by definition of } s_\b 
\end{eqnarray*}
as required.

\nl (iii). The equations are necessary as they appeared under 
(\ref{00eq}).
\end{proof}

\np
The proposition enables us to describe
an algorithm computing the $T_{i,\b}$.

\begin{Alg}\label{alg}
The Hecke algebra elements $T_{i,\b}$ of Theorem \ref{adethmgen2} can be 
computed 
as follows by using Table \ref{algtable}.
\begin{enumerate}[(i)]
\item
If $i\not\in\Supp(\b)$, then, in accordance with (\ref{knotleq}), set $T_{i,\b} = 0$.
\\
{From} now on, assume $i\in \Supp(\b)$.
\item
If ${\rm ht}(\b)\le2$, Equations
(\ref{11eq2}) and (\ref{11eq1}),
that is, the second and third lines of Table \ref{algtable},
determine
$T_{i,\b}$. 

{From} now on, assume 
${\rm ht}(\b)>2$.
We proceed by recursion, expressing $T_{i,\b}$
as an $Z_0^{(0)}$-bilinear combinations of 
$T_{k,\c}$'s with ${\rm ht}(\c)<{\rm ht}(\b)$.
\item
If $(\a_i,\b)=1$, in accordance with Corollary \ref{jastline}(ii), set $T_{i,\b} = m d_{\a_i}^{-1}s_\b^{-1}s_is_\b d_{\b}$.
\\
{From} now on, assume $(\a_i,\b)\in \{0,-1\}$.

\item
Search for a $j\in\{1,\ldots,n\}$ such that
$(\a_i,\a_j) = 0$ and $(\a_j,\b)=1$.
If such a $j$ exists, then
$\b-\a_j\in\Phi$ and (\ref{commlaweq})
expresses $T_{i,\b} $ as a multiple of $T_{i,\b-\a_j}$.
\item
So, suppose there is no such $j$.
There is a $j$ for which 
$\b-\a_j$ is a root, so $(\a_j,\b)=1$.
As $(\a_i,\b)\ne 1$,
we must 
have $i\sim j$.  According as
$(\a_i,\b)=0$ or $-1$, the identities
(\ref{01eq1}) or (\ref{1-1eq1})
express $T_{i,\b}$ 
as a $Z_0^{(0)}$-bilinear combination of $T_{i,\b-\a_j}$ and some $T_{j,\c}$
with ${\rm ht}(\c)<{\rm ht}(\b)$.
\end{enumerate}
This ends the algorithm.
Observe that all lines of Table \ref{algtable} have been used, with 
(\ref{01eq2}) implicitly in (iii).
\end{Alg}

\np The algorithm computes a Hecke algebra element for each $i,\b$ 
based
on Table \ref{algtable}, showing that there is at most one solution 
to
the set of equations. The next result shows that the computed 
Hecke algebra elements do indeed give a solution.

\np
\begin{Prop}\label{uniquesol}
The equations of Table \ref{algtable} have a unique solution.
\end{Prop}

\np
\begin{proof}

We will first show that the Hecke algebra elements $T_{i,\b}$ defined by
Algorithm \ref{alg} are well defined by the algorithm and then that they
satisfy the equations of Table \ref{algtable}.   Both
assertions are proved by induction on $\het(\b)$, the height of $\b$.

If $\b$ has height $1$ or $2$, $T_{i,\b}$ is 
chosen in Step (i) if $\b=\a_j$ with $j\neq i$ and in Step (ii) otherwise.  
Indeed there is a unique solution.

Now assume
$\het(\b)\ge3$.
Suppose first that $T_{i,\b} $ is determined in Step (iii).  This means that 
$(\a_i,\b)=1$. This is 
unique as it is a closed form.

We now suppose that $T_{i,\b}$ is chosen in Step (iv).  This means 
there is a $j$ for which $(\a_{i},\a_{j})=0$ and $(\a_j,\b)=1$. 
We must show that if there are two such $j$ the result is the same.  
Suppose there are
distinct $j$ and $j'$ for which $(\a_j,\b)=(\a_{j'},\b)=1$ and
$(\a_j,\a_i)=(\a_{j'},\a_i)=0$. Then by our definition
$T_{i,\b}=h_{\a_i,j'}^{-1}T_{i,\b-\a_{j'}}$ and we must show that
$T_{i,\b}=h_{\a_i,j}^{-1}T_{i,\b-\a_{j}}$. If $j\sim j'$, then
$(\b-\a_j,\a_{j'})=2$ and $\b=\a_j+\a_{j'}$ has height $2$. This means 
we can
assume $j\not\sim j'$. Then $(\b-\a_j,\a_{j'})=1$ and
$(\b-\a_{j'},\a_j)=1$. In particular, $\b-\a_j-\a_{j'}$ is also a 
root.
Now apply (\ref{commlaweq}) and the induction hypothesis to see
$T_{i,\b-\a_j}=h_{\a_i,j'}^{-1} T_{i,\b-\a_j-\a_{j'}}$ and
$T_{i,\b-\a_{j'}}=h_{\a_i,j}^{-1}T_{i,\b-\a_j-\a_{j'}}$, and so
by (\ref{H-com-rel}), we find
$h_{\a_i,j}^{-1}T_{i,\b-\a_j}=h_{\a_i,j'}^{-1}T_{i,\b-\a_{j'}}$. This shows 
the definitions are the same with either choice.  

We may now assume that $T_{i,\b}$ was chosen in Step (v). If $j$ 
is the one chosen in Step (v), then $T_{i,\b}$ was chosen to 
satisfy (\ref{01eq1}) or (\ref{1-1eq1}). 
Suppose now that there is another 
index $j'$ which was used in Step (v) to define $T_{i,\b}$. For these 
the conditions are $(\a_j,\b)=(\a_{j'},\b)=1$ and 
$(\a_i,\a_j)=(\a_i,\a_{j'})=-1$. Clearly $j\not\sim j'$ for 
otherwise there would be a triangle in the Dynkin diagram $M$. 
Therefore, $(\a_{j'},\b-\a_j) = 1$, and so $\b-\a_j-\a_{j'}$ is a root.
We distinguish according to the two possiblities for $(\a_i,\b)$.

Assume first $(\a_i,\b) = 0$. Then, $(\a_{i},\b-\a_j-\a_{j'})
= 2$, and so $\b = \a_i+\a_j+\a_{j'}$.  By using (\ref{01eq1}),
with either $j$ or with $j'$, we find $T_{i,\b} = m^2$, independent of the choice of
$j$ or $j'$.

Next assume $(\a_i,\b) = -1$.  Then
$(\a_i,\b-\a_j-\a_{j'}) = 1$, so
$\c = \b-\a_j-\a_{j'}-\a_i$ is a root.
We need to establish that
the result of application of (\ref{1-1eq1}) 
to $T_{i,\b}$ does not depend on the choice $j$ or $j'$.
We do so by showing that the result can be expressed in an expression
symmetric in $j$ and $j'$. Observe that $\c$ is an expression symmetric in $j$
and $j'$.
The expression of $T_{i,\b}$ obtained by applying (\ref{1-1eq1}) to $j$ is
\begin{eqnarray}\label{jlaatsteq}
 && T_{j,\b-\a_j}h_{\b-\a_j,i}+ m T_{i,\b-\a_j}.
\end{eqnarray}
By
(\ref{01eq1}),
the second summand of the right hand side
equals
\begin{eqnarray*}
mT_{i,\b-\a_j} & = & mT_{j',\c}+m^2T_{i,\b-\a_j-\a_{j'}}.
\end{eqnarray*}
For the first summand of (\ref{jlaatsteq}) we find
\begin{eqnarray*}
T_{j,\b-\a_j}h_{\b-\a_j,i} &=&
 h_{\a_j,j'}^{-1}T_{j,\b-\a_j-\a_{j'}}h_{\b-\a_j,i}
 \qquad \hfill \phantom{mxucchmucccchmorre}
        \mbox{ by (\ref{commlaweq}) } 
\\
 &=&
 h_{\a_j,j'}^{-1}
 \left(T_{i,\c}h_{\c,j}+ m T_{j,\c}\right)
 h_{\b-\a_j,i}.
 \qquad \hfill \phantom{orre}
\mbox{ \ \ \ \ \ \ \ \ \ \ \ \ \ \  by (\ref{1-1eq1}) } 
\end{eqnarray*}

Expanding (\ref{jlaatsteq}) with these expressions, we find
by use of $h_{\a_j,j'}= h_{\a_{j'},j}$ (see (\ref{H-0-1-rel}),
$h_{\c,j'} = h_{\b-\a_j,i}$ (see (\ref{H-dist2-rel})), and (\ref{commTH}),
\begin{eqnarray*}
&& h_{\a_j,j'}^{-1} T_{i,\c}h_{\c,j}h_{\b-\a_j,i}
 + m \left(h_{\a_j,j'}^{-1}T_{j,\c}h_{\b-\a_j,i} + T_{j',\c}\right)
 + m^2 T_{i,\b-\a_j-\a_{j'}}=\\
&&h_{\a_j+\a_{j'}+\a_i,i}^{-1} T_{i,\c}h_{\c,j}h_{\c,j'}
 + m \left(T_{j,\c} + T_{j',\c}\right)
 + m^2 T_{i,\b-\a_j-\a_{j'}}.
\end{eqnarray*}
Since $h_{\c,j}$ and $h_{\c,j'}$ commute, cf.\ (\ref{H-com-rel}), the result is
indeed symmetric in $j$ and $j'$.
This shows that the algorithm gives unique Hecke algebra elements $T_{i,\b}$.

We now show that the relations of Table \ref{algtable} all hold for $T_{i,\b}$
as computed by the algorithm.  If the height of $\b$ is one or two the values
are given by (\ref{knotleq}) and (\ref{11eq2}) of the table and none of the other
relations holds as there are no applicable $j$.

We consider each of the remaining relations, one at a time, and show that 
each holds by assuming the relations all hold for roots of lower height.

If $(\a_i,\b)=1$ the value of $T_{i,\b}$ is given in Step (iii).  The relevant equations are (\ref{commlaweq}) 
and (\ref{01eq2}).  The proof of Corollary \ref{jastline}(ii)
shows that both equations are satisfied by the closed formula which is the
outcome of our algorithm.

We have yet to check (\ref{01eq1}) and (\ref{1-1eq1}) in which case $(\b,\a_i) $ is $0$ or $-1$.  Notice 
(\ref{commlaweq}) and (\ref{01eq2}) require $(\a_i,\b)=1$ and do not apply here.  In 
these cases $T_{i,\b}$ is chosen in Step (iv) or Step (v).     

Suppose 
first $T_{i,\b} $ was chosen by Step (iv). In this case there is a $j'$ with
$(\a_{j'},\b)=1$, $(\a_i,\a_{j'})=-1$.
As $T_{i,\b}$ is determined by Step (iv) of the algorithm, 
$T_{i,\b}=h_{\a_i,j'}^{-1}T_{i,\b-\a_{j'}}$.  We have already seen that 
this is independent of the choice of $j'$ and so if there is another 
$j$ for which $(\a_j,\b)=1$ with $(\a_i,\a_j)=1$, (\ref{commlaweq}) holds.  To check 
(\ref{01eq1}) we suppose there is a $j$ for which $(\a_i,\b)=1$ with $(\a_i,\a_j)=-1$.  We
must have $j\not\sim j'$, for otherwise we would again
be in the height $2$ case. In order to obtain
(\ref{01eq1}) we must show that 
$$
h_{\a_i,j'}^{-1}T_{i,\b-\a_{j'}}=T_{j,\b-\a_i-\a_j}+m T_{i,\b-\a_j}.
$$
As for the left hand side, $(\b-\a_{j'},\a_j)=1$ and $(\a_i,\a_j)=-1$, 
so by
(\ref{01eq1}), we have
$$
h_{\a_i,j'}^{-1}T_{i,\b-\a_{j'}}=h_{\a_i,j'}^{-1}T_{j,\b-\a_{j'}-\a_j-\a_i}
+mh_{\a_i,j'}^{-1}T_{i,\b-\a_{j'}-\a_j}.
$$
As for the right hand side, as $(\a_j,\a_{j'})=0$, we can
use (\ref{commlaweq}) to obtain
$T_{j,\b-\a_i-\a_j}=h_{\a_j,j'}^{-1}T_{j,\b-\a_j-\a_i-\a_{j'}}$ and
$T_{i,\b-\a_j}=h_{\a_i,j'}^{-1}T_{i,\b-\a_j-\a_{j'}}$, and so the right hand side 
equals the left hand side if $ h_{\a_j,j'}= h_{\a_i,j'} $.
But this is
(\ref{H-jjprime-rel}).

We have yet to consider the case $(\a_i,\b)=-1$, when $T_{i,\b}$ is chosen in Step (iv).  Suppose 
$j'$ is the choice used in Step (iv).    
As we saw in the case $(\a_i,\b)=0$, (\ref{commlaweq}) holds for any $j$ with $(\a_j,\b)=1$ and with 
$(\a_i,\a_j)=0$ by the uniqueness of the definition of $T_{i,\b}$.  
We need to treat the case $(\a_j,\b)=1$ with $(\a_i,\a_j)=-1$ and show (\ref{1-1eq1}) 
holds.  In particular we need to show 
$$
 h_{\a_i,j'}^{-1}T_{i,\b-\a_{j'}}= T_{j,\b-\a_j} h_{\b- \a_j,i}+mT_{i,\b-\a_j}.
$$
Use (\ref{1-1eq1}) on the left hand side to get 
$$h_{\a_i,j'}^{-1}T_{j,\b-\a_{j'}-\a_j}h_{\b-\a_j-\a_{j'},i}+mh_{\a_i,j'}^{-1}T_{i,\b-\a_j-\a_{j'}}.$$
On the right hand side use (\ref{commlaweq}) to get 
$$h_{\a_j,{j'}}^{-1}T_{j,\b-\a_j-\a_{j'}}h_{\b- \a_j,i}+mh_{\a_i,j'}^{-1}T_{i,\b-\a_j-\a_j}.$$
The needed equation will hold provided $h_{\a_i,j'}=h_{\a_j,j'}$ and $h_{\b-\a_j-\a_{j'},i}=h_{\b-\a_j,i}$.  
The first is (\ref{H-jjprime-rel}) and the second is (\ref{H-orth-rel}).

This shows that all the equations are satisfied if $T_{i,\b} $ is chosen 
in Step (iv).   But if $T_{i,\b}$ was chosen in Step (v) we have already checked any two choices 
of $j$ give the same answer for (\ref{1-1eq1}) and so this equation is satisfied also.  We have now
shown all the relations in Table \ref{algtable} hold.    
\end{proof}

\np At this point we have established the existence of a linear representation
$\s$ of $A$ on $V^{(0)}$.
We need some properties of projections which have already arisen in
\cite{CW}.  In particular let $f_i = ml^{-1}e_i$.  The following lemma shows these
elements are multiples of projections.

\begin{Lm}\label{eiproj}
The endomorphisms $\sigma(f_i)$ of $V^{(0)}$ satisfy
$$
\sigma(f_i)x_\b = 
\left\{
\begin{tabular}{rcl}
$(l^{-2}+ml^{-1}-1)x_{\a_i}$  && if $(\a_i,\b)=2 $, \\
$l^{-1}x_{\a_i}T_{i,\b}(h_{\b,i}+m+l^{-1})$ && if $(\a_i,\b)=0$,\\
$l^{-1}x_{\a_i}(T_{i,\b+\a_i}+l^{-1}T_{i,\b})$ && if $(\a_i,\b)=-1 $,\\
$l^{-1}x_{\a_i}(T_{i,\b-\a_i}+(m+l^{-1})T_{i,\b}) $ && if $(\a_i,\b)=1 $.
\end{tabular}
\right.
$$
In particular, $\sigma(f_i)x_\b\in x_{\a_i}l^{-1}Z_0^{(1)}[l^{-1}]$ if $\b\ne \a_i$
and
$\sigma(f_i)x_{\a_i} \in x_{\a_i}(-1+ l^{-1}Z_0^{(1)}[l^{-1}])$.     
\end{Lm}

\begin{proof}
Suppose first $(\a_i,\b)=2$ in which case $\b=\a_i$.  
Using the definition of $\sigma$ and (\ref{11eq2}) gives $\sigma_ix_{\a_i}=l^{-1}x_{\a_i}$.  
Now $\sigma(f_i)x_{\a_i}=(l^{-2}+ml^{-1}-1)x_{\a_i}$.  

Suppose $(\a_i,\b)=0$.  Then $\sigma_ix_\b=x_\b h_{\b,i}+l^{-1}x_{\a_i}T_{i,\b}$.  Now 
$\sigma_i^2x_\b=x_\b h_{\b,i}^2 +l^{-1}x_{\a_i}T_{i,\b}h_{\b,i}+l^{-2}x_{\a_i}T_{i,\b}$.  
Evaluating $\sigma(f_i)$ on $x_{\a_i}$ and using the Hecke algebra quadratic relation for $h_{\b,i}$ gives that 
the coefficient of $x_\b$ is 
$0$.  Adding the other terms gives $l^{-1}x_{\a_i}T_{i,\b}(h_{\b,i}+m+l^{-1})$ as 
stated.   

Suppose $(\a_i,\b)=-1$.  Now $\sigma_ix_\b=x_{\b+\a_i}-mx_\b+l^{-1}x_{\a_i}T_{i,\b}$.  Applying 
$\sigma_i$ again gives $\sigma_i^2x_\b=x_{\b}+l^{-1}x_{\a_i}T_{i,\b+\a_i}-m(x_{\b+\a_i}-mx_\b+l^{-1}x_{\a_i}T_{i,\b})+l^{-2}x_{\a_i}T_{i,\b}$.  Again adding gives the result.  

If $(\a_i,\b)=1$, $\sigma_i x_\b=x_{\b-\a_i}+l^{-1}x_{\a_i}T_{i,\b}$.  Now $\sigma_i^2x_\b=
x_{\b}-mx_{\b-\a_i}+l^{-1}x_{\a_i}T_{i,\b-\a_i}+l^{-2}x_{\a_i}T_{i,\b}$.  Adding and again using 
the quadratic relation gives the 
result.  

The final statement follows from the fact that 
the $T_{i,\gamma}$ and $h_{\b,i}$ belong to $Z_0^{(1)}[l^{-1}]$ (that is, there  is no
$l$ involved).
\end{proof}

\nl
{\bf Proof of Theorem \ref{adethmgen2}.}
In view of Proposition \ref{finalTiRelations} we need only check (D1), (R1), (R2),
and that $\sigma(e_ie_j)=0$ for $i\not \sim j$.  But (D1) is just the
definition.  
By Lemma \ref{eiproj} we know $\sigma(e_i)x_\b$ is in the space spanned by
$x_{\a_i}$.  Now (R1) follows as $\sigma_ix_{\a_i}=l^{-1}x_{\a_i}$.  For
$i\not \sim j$ we know $\sigma(e_ie_j)=\sigma(e_je_i)$.  By Lemma \ref{eiproj}
this is in $x_{\a_i}Z_0^{(0)}$ and also in $x_{\a_j}Z_0^{(0)}$, and so
it is $0$.  As for (R2) again $\sigma(e_i)x_{\b}$ is a multiple of $x_{\a_i}$.
Now $\sigma_jx_{\a_i}=x_{\a_i+\a_j}-mx_{\a_i}$.  Lemma \ref{eiproj} gives
$\sigma(f_i)(x_{\a_i+\a_j}-mx_{\a_i})=
x_{\a_i}(l^{-1}(m+l^{-1})m-(l^{2}+ml^{-1}-1)m=mx_{\a_i}$.  Now scaling to get
$\sigma(e_i)$ gives the result.  We have shown that Theorem \ref{adethmgen2}
holds.

\np We now show how to construct irreducible representations of $B$ which have
$I_2$ in the kernel.

\begin{Lm} \label{zixa=xaz0-lm}
For each node $i$ of $M$,
we have $\sigma(Z_i^{(0)})x_{\a_i} = x_{\a_i} Z_0^{(0)}$.
\end{Lm}

\begin{proof}
For $j$ and $i$ adjacent nodes, the following computation shows
$\sigma_i\sigma_jx_{\a_i} = x_{\a_j}$.
\begin{eqnarray*} \sigma_i \sigma_j x_{\a_i} &=& \sigma_i(x_{\a_i+\a_j}-mx_{\a_i})
       =x_{\a_j}+l^{-1}T_{i,\a_i+\a_j}x_{\a_i}-ml^{-1}x_{\a_i} \\
                                  &=& x_{\a_j}+l^{-1}x_{\a_i}m-ml^{-1}x_{\a_i} 
                                   =x_{\a_j}.
\end{eqnarray*}
By induction on the length of a path from $i$ to $k$ in $M$, this gives
\begin{eqnarray}\label{xaitoj-eq}
\sigma(\wh{w_{ik}})x_{\a_i} &=&   x_{\a_k}.
\end{eqnarray}

Therefore, for $j$ and $k$ distinct nonadjacent nodes of $M$,
$$
x^{-1} \sigma(\wh{w_{ki}}\wh{j}\wh{w_{ik}}e_i)x_{\a_i} 
= 
 \sigma(\wh{w_{ki}}\wh{j})x_{\a_k}
=
\sigma(\wh{w_{ki}})\sigma_jx_{\a_k} \\
=
 \sigma(\wh{w_{ki}})x_{\a_k}h_{\a_k,j}
= x_{\a_i}h_{\a_k,j} .
$$
As $\sigma(Z_i^{(0)})$ is generated by elements of the form
$ \sigma(\wh{w_{ki}}\wh{j}\wh{w_{ik}}e_i) $, it follows that
$\sigma(Z_i^{(0)})x_{\a_i} \subseteq x_{\a_i} Z_0^{(0)}$.  Note it follows 
from Lemma \ref{eiproj} that $x^{-1}\sigma(e_i)x_{\alpha_i}=x_{\alpha_i}$.

As for the converse, this follows from Lemma \ref{Hbeta-lm}(ii), which implies that
$Z_0^{(0)}$ is generated by $h_{\a_k,i}$, for $i\not\sim k$, $i\ne k$.
(For, by definition, $Z_0^{(0)}$ is generated by $\wh{C}$ mod $I_2$.)
\end{proof}

Suppose $\theta$ is any representation of $Z_0$, acting on a vector space $U$
over $K$, where $K = \Q(r)$, or an algebraic extension thereof.  Then we can form a
representation of $B$ on the vector space $V\otimes_{Z_0}U$ over $K(l)$
which is the direct sum of vector spaces $x_\b U$ where each is a vector space
isomorphic to $U$.  Let $V$ be the representation space of Theorem
\ref{adethmgen2}.  For each $i$ define an action of $\sigma_i$ on
$V\otimes_{Z_0}U$ by letting elements of $Z_0$ act directly on $U$.  In
particular, $\sigma_ix_{\a_i}u=l^{-1}x_{\a_i}u$; if $(\a_i,\b)=0$, then
$\sigma_ix_\b u=x_\b \theta(h_{\b,i})u+l^{-1}x_{\a_i}\theta(T_{i,\b})u$; for
$(\a_i,\b)=1$ we have $\sigma_i(x_\b
u)=x_{\b-\a_i}u+l^{-1}x_{\a_i}\theta(T_{i,\b})u$ and if $(\a_i,\b)=-1$ we have
$\sigma_ix_\b u=x_{\b+\b_i}u-mx_\b u +l^{-1}x_{\a_i}\theta(T_{i,\b})u$.  This
is a representation by Theorem \ref{adethmgen2}.  Denote it $\Gamma_\theta$.

\begin{Lm}  \label{irr-rep-lm}
If $\theta $ is an irreducible representation of $Z_0^{(0)}$, then the representation 
$\Gamma_\theta$ is also irreducible.   
For inequivalent representations $\theta$, $\theta'$, the resulting
representations $\Gamma_\theta$ and
$\Gamma_{\theta'}$ are also inequivalent.
\end{Lm}

\begin{proof}  Suppose $V_1$ is a proper nontrivial invariant subspace of $V\otimes_{Z_0}U$.  We show first 
  that $\sigma(f_i)V_1=0$ for all nodes $i$ of $M$.  By Lemma
  \ref{eiproj},  $\sigma(f_i)V\otimes_{Z_0}U$ is in $x_{\a_i}\theta(Z_0^{(0)})U$ which is
  in $x_{\a_i}U$.  This means that $\sigma(f_i)V_1$ is in $x_{\a_i}U$.  Suppose
there is a node $i$ with  $\sigma(f_i)V_1$ nonzero.  This means there is a nonzero element of $u\in
  U$ such that $x_{\a_i}u\in V_1$.  In Lemma \ref{zixa=xaz0-lm}, we have seen
  that $Z_i^{(0)}x_{\a_i}=x_{\a_i}Z_0^{(0)}$.  Hence $x_{\a_i}\theta(Z_0^{(0)})u
  = Z_i^{(0)}x_{\a_i}\subseteq V_1 $.  But $\theta$ is irreducible and so all
  of $x_{\a_i}U$ is contained in $V_1$.

By Lemma \ref{zixa=xaz0-lm}, $x_{\a_k}U$ is in $V_1$ for
all $k$.  We show by induction on the height of a positive root
$\het(\b)$ that $x_\b U$ is in $V_1$.  Assume $\het(\b)\ge2$. Choose a node $j$
with $\b=r_j(\b-\a_j)$.  By induction, $x_{\b-\a_j}U$ is in $V_1$.  But for
each $u\in U$, the vector $\sigma_jx_{\b-\a_j}u$ is a sum of $x_\b u$ and
vectors already known to be in $V_1$ and so $x_\b U$ is in $V_1$.  But this
means all of $V\otimes_{Z_0}U$ is in $V_1$, contradicting that $V_1$ is proper.
This shows $\sigma(f_i)V_1=0$ for each node $i$.  

As $V_1$ is invariant, its image
$\sigma(\wh{w_{\b,j}}f_j\wh{w_{\b,j}}^{-1})V_1$ under a conjugate of
$\sigma(f_i)$ is also trivial.  We
will derive from this that $V_1$ is $0$.  To this end, choose an order on
$\Phi^+$ that is consistent with height.  For each $\b$ choose a node $j(\b)$ in
the support of $\b$.  Notice that Lemma \ref{eiproj} shows that the image of
$\sigma(f_i)$ is in $x_{\a_i}Z_0^{(0)}$.  Let $L$ be the matrix whose rows and columns are
indexed by $\Phi^+$ in the fixed order and whose $\b,\c$ entry is the
coefficient of $x_{\b}$ in
$\sigma(\wh{w_{\b,j(\b)}}f_{j(\b)}\wh{w_{\b,j(\b)}}^{-1})x_\c$.  This means the
entries are elements of $\theta(Z_0^{(0)})$.  As each
$\sigma(\wh{w_{\b,j(\b)}}f_{j(\b)}\wh{w_{\b,j(\b)}}^{-1})V_1=0$, we have $LV_1=0$.

Observe that $L$ can be viewed as a matrix with entries in $K[l^{-1}]$ by interpreting
the entries from $\theta(Z_0^{(0)})$ as submatrices over $K[l^{-1}]$.
We claim that $L$ is nonsingular.
By the Lawrence-Krammer action rules,
the $\b,\c$ entry of $L$ mod $l^{-1}$ is readily seen to be the coefficient of $x_{\a_{j(\b)}}$ in
$\sigma(f_{j(\b)}\wh{w_{\b,j(\b)}}^{-1})x_\c$. 
If $\b=\c$, then this coefficient is equal to $-1$ 
modulo $l^{-1}$, and if $\b $ is less than $\c$ in the given order, then
there is no summand $x_{\a_{j(\b)}}$ present in the expansion of $\sigma(\wh{w_{\b,j(\b)}}^{-1})x_\c$
and so the $\b,\c$ coefficient of $L$ is $0$.
This means $L$ modulo $l^{-1}$ is lower-triangular with $-1$ on the diagonal, whence non-singular.

Therefore, the equality $LV_1 = 0$ implies $V_1=0$.  We conclude that there is
no invariant subspace and the representation is irreducible.

Finally, we argue that inequivalent $\theta$ lead to inequivalent
$\Gamma_\theta$. To this end we consider the trace of 
each element $\wh{w_{ki}}\wh{z}\wh{w_{ik}}e_i$ of $Z_i$ in $\Gamma_\theta$,
where $z$ is in $W_{k^\perp}$.
By Lemma \ref{eiproj}, the only contributions to the trace occur for
vectors in $x_{\a}\theta(Z_0)$, and, in view of 
Lemma \ref{zixa=xaz0-lm}, this contribution is
$m^{-1}(l^{-1}+m-l^{-1})\tr(\theta(d_{\a_k}^{-1}\wh{z}d_{\a_k}))$.
Since $d_{\a_k}^{-1}\wh{z}d_{\a_k}$, for $k$ a node of $M$ and $z \in
W_{k^\perp}$, span $Z_0$ over $K(l)$, these values uniquely determine $\theta$.
\end{proof}

\np
With these results in hand we are now ready to show that the dimension of $I_1/I_2$ is at least the 
dimension we need for Theorem \ref{main-thm}.

\nl
{\bf Proof of Theorem \ref{main-thm}.}
In Theorem \ref{irr-rep-lm} we have constructed irreducible representations $\Gamma_\theta$
of $B/I_2$
of dimension $|\Phi^+|\dim\theta$ for any irreducible representation $\theta$ of $Z_0$.
Since $I_1$ is not in the kernel of these representations,
they are irreducible representations of $I_1/I_2$.
Moreover, $Z_0$, being a Hecke algebra over $\Q(l,m)$ of spherical type, is semi-simple, so summing the squares of the 
dimensions of the irreducibles of $Z_0$ gives $\dim(Z_0)$.
Hence the dimension of $I_1/I_2$ is at least $ |\Phi^+|^2\dim(Z_0)$. 
By Theorem
\ref{doubleI}, this is also an upper bound for the dimension, whence equality.  The 
semisimplicity follows as $B/I_1$, being the Hecke algebra of type $M$,
is semisimple, and the sum of the squares of the irreducible representations 
of $I_1/I_2$ is the dimension of $I_1/I_2$.  

\np To end this section, we observe that the usual Lawrence-Krammer
representation is the representation $\Gamma_\theta$, where $\theta$ is the
linear character of $Z_0$ determined by $\theta(h_{\b,i}) = r^{-1}$ for all
pairs $(\b,i)\in \Phi^+\times M$ with $(\a_i,\b) = 0$.

\section{Consequences and Conjectures} \label{sec:concon}
This section gives some consequences of the main results of the previous
sections, as well as some of our ideas about the general structure of BMW
algebras.

\subsection{Global structure of BMW algebras}\label{global-subsec}
Indications for the validity of our theorems were first found by experimental
computations in GBNP, \cite{GBNP}.  However, the sheer size of the algebras
involved makes the computations difficult.  For instance, the dimension of
$I_1/I_2$ in $B(\E_8)$ is equal to $41803776000$.

Nevertheless, some experimenting with $B(\D_4)$ and knowledge of the classical
BMW algebra $B(\A_n)$
lead us to conjecture that, if $J$ is a coclique of $M$ of size $i>1$,
then $I_J$ is an ideal properly contained in $I_{i-1}$.

If $J$ and $K$ are conjugate by an element $w\in W$, then
as we have seen in Proposition \ref{Ideal-I-prop}(ii), the ideals $I_J$ and $I_K$ coincide.
Computations in $B$ of type $\D_4$
show that for $J$ and $K$ of size 2 but in distinct orbits,
we find distinct ideals $I_J = Be_JB$, $I_K = Be_KB$.
Also the pattern that, for each
coclique $J$ of size $i$, we have $I_J/I_{i+1}= Be_JB/I_{i+1} = 
\wh{D_J}Z_J\wh{D_J}^{op}/I_{i+1}$ for a suitable set $D_J$ of coset
representatives of the stabilizer of $ \{r_j\mid j\in J\}$ in $W$ and a 
subalgebra $Z_J$ of $B$ isomorphic to a suitable subtype $C_J$ of $M$.
Thus, we expect that $\dim(I_J/I_{i+1}) $ is a multiple of $N^2$ by the order
of a Coxeter group of some subtype $C_J$ of $M$,
where $N$ is the length of the $W$-orbit of $ \{r_j\mid j\in J\}$.
This would imply that the dimension of $B$ be equal to
$$\sum_J N_J^2\,|W(C_J)|.$$
Here $J$ runs over the $W$-equivalence classes of cocliques in $M$, including the empty
set,
with $C_\emptyset = M$ and $N_\emptyset = 1$, so that 
the contribution for $J = \emptyset$ equals $|W|$, the dimension of $B/I_1$,
the Hecke algebra of type $M$.

The conjecture holds for $B(\A_n)$. Here $W$ is known to have
a single orbit on cocliques in $M$ of any given size
$i\in\{1,\ldots,\lceil{n/2}\rceil\}$;
for $J = \{1,3,\ldots,2i-1\}$, the type $C_J$ is the Coxeter type of the
centralizer in $W$
of $\{\a_j\mid j\in J\}$, that is, $C_J = \A_{n-2i}$, and
\begin{eqnarray*}
\dim(I_i/I_{i+1}) &=& 
N_i^2\, (n+1-2i)!\qquad\mbox{ with } \ 
N_i = \binom{n+1}{\underbrace{2,2,\ldots,2}_{i\ \times}}.
\end{eqnarray*}
These formulas also hold for $i=0$ if we write $I_0 = B$ and $N_0= 1$. We
then find
$\dim(B(\A_n)) = \sum_i\dim(I_i/I_{i+1}) =
(2n+1)(2n-1)(2n-3)\cdots 1$, which is known from \cite{Wenzl}.

Our conjecture also holds for $B(\D_4)$.
In $B(\D_4)$, there are three ideals of the form $I_J$ for $J$ of size 2,
namely for $J = \{1,3\}$, $\{1,4\}$, $\{3,4\}$. Each quotient $I_J/I_3$
has dimension $N_J^2\cdot 2$, where $N_J = 6$.
Thus $C_J$ is of type $\A_1$, rather than $\A_1\A_1$, the
parabolic type of the centralizer of two orthogonal roots. This means that a
complication with respect to the type $\A_n$ occurs in that the type $C_J$
is not just the full type of the centralizer of $\{\a_j\mid j\in J\}$ in $W$.
Similarly, $N_{\{1,2,3\}} = 3$, $C_{\{1,2,3\}} = \emptyset$,
and $I_3= I_{\{1,3,4\}}$  has dimension $N_{\{1,2,3\}} ^2\cdot 1 = 9$.
In conclusion, $\dim(B(\D_4)) = |W| + N_1^2 |W(\A_1^3)| + 
3\times N_{\{1,3\}}^2 |W(\A_1)| + N_{\{1,3,4\}}^2 |W(\emptyset)|  = 
192 + 12^2\cdot 8 + 3\cdot 6^2\cdot 2 + 3^2 = 1569$. 

The shrink of $C_J$ for $J$ of size 2 extends to all types $\D_n$ for $n\ge4$.
In $B(\D_n)$ $(n\ge 5)$, there are two conjugacy classes,
one of which has representative $\{n-1,n\}$.
In this case, or rather, in any case where $J$ contains these two end nodes,
the representation of $B(\D_n)$ on $I_J$ factors through a representation
of $B(\A_{n-1})$.
We prove this as
follows.
To begin, we can take $J = \{n-1,n\}$.
We claim that $g_{n}$ acts precisely as $g_{n-1}$. First of all $g_ne_J =
l^{-1} e_J = g_{n-1} e_J$. We proceed to show
$g_n\wh ue_J = g_{n-1}\wh ue_J$ by induction on the length of $u\in W_{\{1,\ldots,n-1\}}$.
Without loss of generality, we may assume $u\in D_{n^\perp,n^\perp}$
(observe that $n^\perp\cap \{1,\ldots,n-1\} = J^\perp\cup\{n-1\}$ in this case, 
so $\wh{g_naub}e_J =\wh a \wh {g_n} \wh u e_J\wh b$ for $a,b\in n^\perp\cap\{1,\ldots,n\}$).
But then, by known properties of the Coxeter group, we have either
$\wh u=g_{n-2}$ or
$\wh u = g_{n-2}g_{n-1}g_{n-3}g_{n-2}$.
As all indices are in $\{n-3,\ldots,n\}$, the identity 
$g_n\wh ue_J = g_{n-1}\wh ue_J$ can be verified in $B(\D_4)$
(after specialization to $n = 4$), where it is easily seen to hold.
So in all cases, $g_n$ acts exactly like $g_{n-1}$, proving that the $B(\D_n)$
representation on $I_J$ factors through the quotient obtained by identifying $g_n$ and
$g_{n-1}$, and so through a BMW algebra of type $B(\A_{n-1})$.
On the basis of observations like these, we conjecture that
the dimension of $B(\D_n)$ is equal to $(2^n+1)(2n-1)!! - (2^{n-1}+1)n!$

\subsection{Parabolic subalgebras and restrictions}
Let $J$ be a set of nodes of $M$.  We will discuss $B_J$, the subalgebra of
$B$ generated by all $g_j$ with $j\in J$.  Clearly, there is a surjective
homomorphism from $B(J)$, the BMW algebra of type $M|_J$ onto $B_J$. We
conjecture however, at least for $M$ of spherical type, that this map is
an isomorphism. It is an easy consequence of Theorem \ref{main-thm} that
this assertion holds modulo $I_2$, in the sense that $B_J/(I_2\cap B_J)$ is
isomorphic to the quotient of $B(J)$ by its ideal $I_2$.

The restriction of the generalized Lawrence-Krammer representation for $B$ on $V$ over $Z_0$ to
$B_J$ is easy to analyze.  For $a:\M\setminus J \to \N$, put
$\Phi^+_{J,a}=\{\b\in\Phi^+\mid C_{\b,k} = a_k \mbox{ for } k\in M\setminus
J\}$ and let $V_{J,a}$ be the subspace of $V$ generated by $x_\b$ with
$\b\in\Phi^+_{J,a}$.  Then it is easily seen from the Lawrence-Krammer action
rules that $V_{J,0}$ is a $B_J$-invariant subspace of $V$, which is isomorphic
to the Lawrence-Krammer representation of $B(J)$, up to an extension of
scalars.  Moreover, the subspace $V_{J,a}+V_{J,0}$ is $B_J$-invariant for any
choice of $a$.  In view of Lemma \ref{eiproj}, the action of $B_J$ on the
quotient $(V_{J,a}+V_{J,0})/V_{J,0}$ factors through the Hecke algebra
$B_J/(I_1\cap B_J)$. We
expect that the particular representations for $B_J$ on $(V_{J,a}+V_{J,0})/V_{J,0}$ can be found by
combinatorics of the root system, similar to the case of type $\A_n$,
discussed in \cite{Wenzl}.

To see how this works in a specific example we consider $B(\D_n)$ with $n\geq
5$ and $J = \{2,3,\dots,n\}$, so we will consider the action of $B_J$
on $V_{J,i}$ for $i=0,1$. Here $\Phi_{J,0}^+$ is the
set of roots $\varepsilon_i \pm \varepsilon_j$ for $2\leq i \leq j \leq n$ and
$\Phi_{J,1}^+$ is the set of roots $\varepsilon_1\pm \varepsilon_j$ for $2\leq
j\leq n$, where $(\varepsilon_i)_{1\le i\le n}$ is an orthonormal basis of
Euclidean $n$-space.  Notice $B_J$ maps the span of $\{x_\beta \mid \b \in
\Phi_{J,0}^+\}$, which is $V_{J,0}$, to itself by the construction for
$B(J)\cong B(\D_{n-1})$.  Also the Hecke algebra $Z_0$ for $B(\D_{n-1})$, which is
$\langle g_2\rangle \times \langle g_4,\dots,g_n\rangle$, can be embedded into
the Hecke algebra $Z_0$ for $B(\D_n)$, which is $\langle g_1\rangle \times
\langle g_3,g_4,\dots,g_n\rangle$, by mapping $g_2$ to $g_1$ and fixing
$\langle g_4,\dots ,g_n\rangle$.  
Furthermore, if $\theta_{res}$ is $\theta $
restricted to $Z_0$ for $B_J$ with this embedding, the resulting
representation of $B(\D_{n-1})$ is $\Gamma_{\theta_{res}}$.  
As mentioned above,
the action of $B(\D_{n-1})$
on the quotient vector space $(V_{J,1}+V_{J,0})/V_{J,0}$
 factors through the Hecke algebra of type
$\D_{n-1}$.  
The represententation then breaks into these two actions with the
action on the quotient being a Hecke algebra action.  
The span $V_{J,1}$ of
the $x_\beta $ for $\beta \in \Phi_{J,1}^+$, is not invariant but using
semisimplicity there is an invariant subspace giving this representation.
This gives a branching rule from $B(\D_{n})$ to $B(\D_{n-1})$.

\subsection{The Brauer algebra}
Let $E$ be the subring $\Q(x)[l^\pm]$ of $\Q(l,x)$.  We conjecture that there is a
subalgebra $B^{(0)}$ of $B$ defined over $E$ containing a spanning set of $B$
with the property that after transition modulo $(l-1)$ we obtain a monomial
algebra whose basis can be described in terms of the root system of type $M$.
For $B$ of type $\A_n$ it is the well-known Brauer algebra, introduced in
\cite{brauer}.  We expect the conjectured basis $\bigcup_J
\wh{D_J}\wh{W_{C_J}}\wh{D_J}^{op}$ of $B$ discussed in \S\ref{global-subsec},
to be a monomial basis mod $E$ for the Brauer algebra.  Its
elements should correspond to pictures, which consist of triples consisting of
two sets of
orthogonal roots, both $W$-conjugate to $\{\a_j\mid j\in J\}$, and an element of
$W(C_J)$, a Coxeter group in a quotient of the centralizer of $J$ in $W$.
This correspondence is well known for type $\A_n$.  The basis of $I_1/I_2$
found in Theorem \ref{main-thm} can be used to establish the validity of this
conjecture for $B/I_2$.

\subsection{Conclusion}
For Coxeter diagrams that are not simply laced, we expect a natural BMW
algebra to exist as well.  For type ${\rm B}_n$, an approach is given in
\cite{haring}. More generally, by means of a folding $\phi: M\to M'$ of Coxeter diagrams, a
BMW algebra of spherical type $M'$ could be constructed as the subalgebra of
$B(M)$ generated by suitable products of $g_i$ for $g_i\in \phi^{-1}(a)$, one
for each $a\in M'$, in much the same way the Artin group of type $M'$ in
embedded into the one of type $M$, see \cite{crisp}. However, further research
is needed to see if this definition is independent (up to isomorphism) of the
choice of $\phi$ for fixed $M'$, as well as to find an intrinsic definition of
this algebra.

The BMW algebras of type $\A_n$ play a role in algebraic topology, in
particular, in the theory of knots. The versions of spherical type $\ADE$ are
related to the topology of the quotient space by $W$ of the complement of the
union of all reflection hyperplanes in the complexified space of the
reflection representation of $(W,R)$.  After all, by \cite{Brieskorn},
the Artin group $A$ is the fundamental group of this space. A direct
relationship, for instance, a definition of the BMW algebra in terms of this
topology, would be of interest.

Brauer algebras play a role in tensor categories for the representations of
classical Lie groups, and the corresponding BMW algebras seem to play
a similar role for the related quantum groups. It is conceivable that
the new BMW algebras constructed here play a similar role for the tensor
categories of representations of quantum groups for the other types.

\end{document}